\documentclass[preprint,12pt]{elsarticle}

\usepackage{graphicx,amsmath,amsfonts,comment}

\def\bfm#1{\boldsymbol{#1}} 
\def\RR{\mathbb{R}}
\def\NN{\mathbb{N}}

\def\Re{\operatorname{Re}}
\def\Im{\operatorname{Im}}
\def\conj#1{\overline{#1}}
\newcommand{\compI}{\operatorname{i}} 

\newtheorem{thm}{Theorem}
\newtheorem{lem}{Lemma}

\newtheorem{rem}{Remark}

\newproof{pf}{Proof}
\newproof{pot}{Proof of Theorem \ref{thm2}}
\journal{Computer Aided Geometric Design}

\biboptions{authoryear}

\begin{document}

\begin{frontmatter}

\title{Arc length preserving approximation of circular arcs by
Pythagorean-hodograph curves of degree seven}

\author[address1,address2]{Emil \v{Z}agar\corref{corauth}}
\ead{emil.zagar@fmf.uni-lj.si}
\cortext[corauth]{Corresponding author}

\address[address1]{Faculty of Mathematics and Physics, University of Ljubljana, 
Jadranska 19, Ljubljana, Slovenia}
\address[address2]{Institute of Mathematics, Physics and Mechanics, Jadranska 19, Ljubljana, Slovenia}

\begin{abstract}
  In this paper interpolation of two planar points,
  corresponding tangent directions and curvatures with 
  Pythagorean-hodograph (PH) curves of degree seven preserving an arc length
  is considered. A general approach using complex representation of
  PH curves is presented and a detailed analysis of the problem for data
  arising from a circular arc is provided. In the case of several 
  solutions some criteria for the selection of the most appropriate one are 
  described and an asymptotic analysis is given. Several numerical examples are included
  which confirm theoretical
  results.
\end{abstract}


\end{frontmatter}


\section{Introduction}\label{sec:introduction}
Interpolation of local planar geometric data, such as points, tangent directions and curvatures, by parametric polynomial
curves is a standard problem in Computer Aided Geometric Design (CAGD) and a common way to
construct parametric objects from given discrete data. If such interpolants are joined together, they form
geometrically continuous splines of order $k$ (or $G^k$ continuous splines),
where $k$ depends on the type of the interpolated data ($k=0$ if only positions of points are
given, $k=1$ if in addition also tangent directions are provided, etc.).
For a detailed survey of geometric interpolation methods
the reader is referred to \cite{HoschekLasser} or \cite{FarinHoschekKim-02-Handbook}. 
However, there are not many results concerning also the
interpolation of some global geometric data, such as the arc length. This becomes extremely important
when some global shape control is needed or methods relying on optimization of curve energies,
such as bending energy, are to be developed. It turned out that there is a specific class of curves which
are of great help to solve such kind of problems, the so called polynomial Pythagorean-hodograph (PH) curves
introduced in \cite{Farouki-Sakkalis-PH-first-90} and comprehensively described in \cite{Farouki-PH-book-2008}.
They are unique among polynomial curves possessing polynomial arc length function. This implies
several nice properties which will be explained in detail later. Some recent results
concerning interpolation of $G^1$ data by PH quintic curves preserving an arc length are in 
\cite{Farouki-quintic-arc-length-16} and they confirm the advantage of PH curves if interpolation
of global geometric data is considered. The author studied the interpolation of two points together
with the corresponding tangent directions and a prescribed arc length. A detailed analysis of the
interpolation problem was done and a simple algorithm relaying basically only on the solution
of the quadratic equation was described. A special type of $G^2$ data interpolation by PH curves 
was considered also in \cite{Farouki-corners-2014}, but the specification of an arc length was not considered. 
In this paper we intend to extend recent results 
to the interpolation of $G^2$ data and an arc length by PH curves of degree $7$. As already guessed 
in \cite{Farouki-quintic-arc-length-16}, it is hard to believe that this problem will posses such simple
solution as in the quintic case. 
We will confirm this prediction after specifying a general problem and later 
concentrate only on interpolation of particular type of data, i.e. the one arising from a circular arc
which will be addressed as circular arc data. It should be noted that even in this case there is no
relevant literature available. Most of the approximation techniques namely consider
the interpolation of local geometric data only 
and use the remaining parameters to minimize the
distance between the interpolant and the circular arc, to minimize the deviation of the curvature, etc.
The results of this type can be found in
\cite{Dokken-Daehlen-Lyche-Morken-90-CAGD}, 
\cite{Goldapp-91-CAGD-circle-cubic}, 
\cite{Lyche-Morken-94-Metric},
\cite{Morken-93-Parametric-Parabola},
\cite{Ahn-Kim-arcs-1997},
\cite{Ahn-Kim-quartic-2007},
\cite{Jaklic-Kozak-Krajnc-Zagar-circle-like-07},
\cite{Jaklic-Kozak-Krajnc-Vitrih-Zagar-conics-13},
\cite{Kovac-Zagar-curvature16},
\cite{Jaklic-circle-CAGD-2016},
\cite{Jaklic-Kozak-2018-best-circle},
\cite{Knez-Zagar-2018-circular-arcs-maxmial-smoothness},
\cite{Vavpetic-Zagar-2018-optimal-circle-arcs},
\cite{Ahn-2019-extreme-points-arcs},
\cite{Vavpetic-2020-optimal-circular-arcs}
and
\cite{Vavpetic-Zagar-21-circle-arcs-Hausdorff},
if we mention just the most important and recent ones. 
Although the proposed algorithms provide good
approximations of circular arcs if the Hausdorff distance is considered as a measure of the error, they 
do not include an arc length in interpolation data. Our approach intend to fill this gap and provide 
interpolants which possess required arc length and remain small Hausdorff distance.

The paper is organized as follows. In Section \ref{sec:preliminaries} some basic properties of
complex representation of PH curves is given. Special attention is given to PH curves of degree $7$
which are  presented in detail and
all quantities needed for a solution of the interpolation problem are derived. In the next section the 
arc length preserving interpolation of $G^2$ data is considered and the system of nonlinear equations
for the general case is presented. Two examples of such interpolation are given showing that the solution
of the problem highly depends on prescribed data. In Section \ref{sec:circular_arc} the
interpolation of circular arc data in canonical position is considered. The system of nonlinear equations
derived for the general case is simplified and a new, simpler system of two nonlinear equations
involving the angle $\alpha$ arising from the circular arc as a parameter are provided. A detailed analysis of the
solvability is done. Two cases are considered. For the first the absence
of real solutions is confirmed,
and for the second one the existence of two real solutions is first proved 
for any $\alpha\in(0,\pi/2]$.
A simple (numerical) procedure to check the existence of four admissible
solutions is described which works well for several particular practically important values  $\alpha$ tested.
In Section \ref{sec:asymptotic_analysis} the existence of four solutions in confirmed for any $\alpha$ small enough.
and asymptotic expansions of solutions are provided. In Section \ref{sec:solution_selection} some criteria
for the selection of the most appropriate solution are described and in the next section several numerical
examples together with the numerical confirmation of the approximation order are given. The paper is concluded
by Section \ref{sec:closure}.

\section{Preliminaries}\label{sec:preliminaries}
Planar PH curves 
are an important subclass of planar parametric polynomial curves. A regular planar parametric
polynomial curve  $\bfm{p}:[0,1]\to\RR^2$ is a PH curve if $\|\bfm{p}'\|$ is a 
polynomial, where $\|\cdot\|$ denotes the standard Euclidean norm on $\RR^2$. This characterization 
implies several important geometric properties of PH curves (\cite{Farouki-PH-book-2008}): 
rational unit tangent, normal, curvature and offset,\dots 
Furthermore, the arc length function of a PH curve is polynomial. 
All these properties make them useful for interpolation
of local geometric data as well as for the interpolation of global geometric quantities, such as an arc length.
Let $\bfm{p}=(x,y)^T$ be a PH curve of degree $n$ where $x'$ and $y'$ are relatively prime polynomials. 
It is known (\cite{Kubota-72}) that polynomials $x'$ and $y'$ can be
expressed in terms of two polynomials $u$ and $v$ as $x'=u^2-v^2$ and $y'=2\,u\,v$
which implies $x'^2+y'^2=(u^2+v^2)^2$ and consequently the arc length function is a polynomial
$u^2+v^2$.
It is often better to use a complex representation of PH curves (\cite{Farouki-complex-PH-94}).
If $n=2\,m+1$ and
\begin{equation}\label{eq:preimage}
  \bfm{w}=\sum_{k=0}^m B_k^m \bfm{w}_k, \quad \bfm{w}_k=u_k+\compI v_k,\ k=0,1,\dots,m,
\end{equation}
where $B_k^m(t)=\binom{m}{k}t^k(1-t)^{m-k}$, $k=0,1,\dots,m$,  are Bernstein basis polynomials, 
then the integral of $\bfm{p}'=\bfm{w}^2$ is a PH curve $\bfm{p}$ of degree $n$.
Furthermore, its unit tangent vector $\bfm{g}$, the curvature $\kappa$  and the arc length $s$
are given by
\begin{equation}\label{eq:properties_w}
  \bfm{g}(t)=\frac{\bfm{w}^2(t)}{\sigma(t)},\quad 
  \kappa(t)=2\frac{\Im\left(\conj{\bfm{w}}(t)\bfm{w}'(t)\right)}{\sigma^2(t)},\quad
  s(t)=\int_0^t\sigma(\tau){\rm d}\tau,
\end{equation}
where $\sigma=\left|\bfm{w}\right|^2$.
Since we will consider PH curves of degree $7$, we shall start by a complex cubic polynomial
$\bfm{w}$ with complex Bernstein coefficients $\bfm{w}_k=u_k+\compI\,v_k$, $k=0,1,2,3$. Integration of
its square results in a PH curve $\bfm{p}$ of degree $7$, which can be written in Bernstein-B\'ezier
form as
\begin{equation*}
  \bfm{p}(t)=\sum_{k=0}^7 B_k^7(t)\,\bfm{p}_k,
\end{equation*}
where
\begin{align}
  \bfm{p}_1&=\bfm{p}_0+\frac{1}{7}\bfm{w}_0^2,\quad 
  \bfm{p}_2=\bfm{p}_1+\frac{1}{7}\bfm{w}_0\bfm{w}_1,\quad
  \bfm{p}_3=\bfm{p}_2+\frac{1}{7}\frac{3\bfm{w}_1^2+2\bfm{w}_0\bfm{w}_2}{5},\nonumber\\
  \bfm{p}_4&=\bfm{p}_3+\frac{1}{7}\frac{9\bfm{w}_1\bfm{w}_2+\bfm{w}_0\bfm{w}_3}{10},\label{eq:p_coeff}\\
  \bfm{p}_5&=\bfm{p}_4+\frac{1}{7}\frac{3\bfm{w}_2^2+2\bfm{w}_1\bfm{w}_3}{5},\quad
  \bfm{p}_6=\bfm{p}_5+\frac{1}{7}\bfm{w}_2\bfm{w}_3,\quad
  \bfm{p}_7=\bfm{p}_6+\frac{1}{7}\bfm{w}_3^2,\nonumber
\end{align}
and $\bfm{p}_0$ is a free complex integration constant. By \eqref{eq:preimage} and
\eqref{eq:properties_w} we obviously have
\begin{equation*}
  \bfm{g}(0)=\left(\frac{\bfm{w}_0}{\left|\bfm{w}_0\right|}\right)^2,\quad
  \bfm{g}(1)=\left(\frac{\bfm{w}_3}{\left|\bfm{w}_3\right|}\right)^2,
\end{equation*}
and by using some basic properties of $\bfm{w}$ we also get 
\begin{equation*}
  \kappa(0)=6\Im\left(\frac{\conj{\bfm{w}}_0\bfm{w}_1}{\left|\bfm{w}_0\right|^4}\right),\quad
  \kappa(1)=-6\Im\left(\frac{\conj{\bfm{w}}_3\bfm{w}_2}{\left|\bfm{w}_3\right|^4}\right).
\end{equation*}
Furthermore, the total arc length $L$ of $\bfm{p}$ is 
\begin{align*}
  L&=\int_{0}^1\left|\bfm{w}(\tau)\right|^2\,{\rm d}\tau
  =\frac{1}{7}\left(\left|\bfm{w}_0\right|^2+\Re\left(\bfm{w}_0\conj{\bfm{w}}_1\right)
  +\frac{2}{5}\Re\left(\bfm{w}_0\conj{\bfm{w}}_2\right)
  +\frac{1}{10}\Re\left(\bfm{w}_0\conj{\bfm{w}}_3\right)
  +\frac{3}{5}\left|\bfm{w}_1\right|^2\right.\\
  &\left.+\frac{9}{10}\Re\left(\bfm{w}_1\conj{\bfm{w}}_2\right)
  +\frac{2}{5}\Re\left(\bfm{w}_1\conj{\bfm{w}}_3\right)
  +\frac{3}{5}\left|\bfm{w}_2\right|^2+\Re\left(\bfm{w}_0\conj{\bfm{w}}_1\right)
  +\left|\bfm{w}_3\right|^2\right).
\end{align*}
These results will now be used in the following section where an interpolation problem of general $G^2$ data
by PH curves of degree $7$ with a prescribed arc length will be considered. 

\section{Arc length preserving interpolation of $G^2$ data}\label{sec:interpolation}

\begin{figure}[htb]
\centering
\includegraphics[width=1\textwidth]{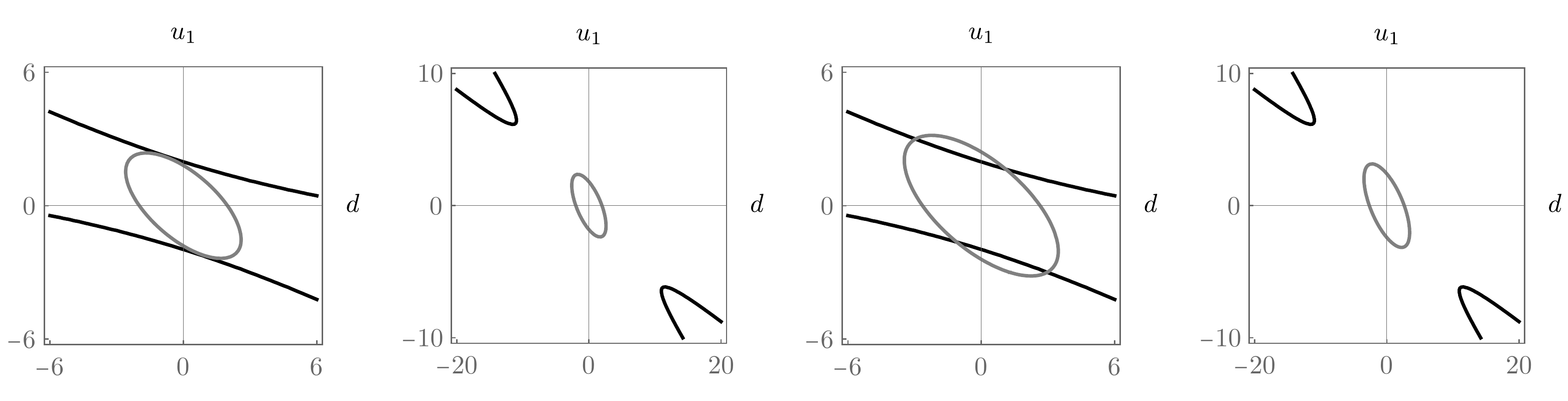}
\caption{Ellipses (gray) and hyperbolas (black) arising from the 
first example with $L=9/8$ (the first and the second figure from the left)
and ellipses and hyperbolas arising from the second
example with $L=2$ (the third and the last figure from the left).}
\label{fig:example12}
\end{figure}

A construction of parametric polynomial curves is usually based on interpolation of particular
geometric data arising from practical observations, such as point positions, tangent directions,
curvatures, etc. As it was already mentioned in the previous section, in some problems also the 
interpolation of global geometric data, such as an arc length, is required. We shall follow the approach in
\cite{Farouki-quintic-arc-length-16}, where the author considered the problem of $G^1$ data 
interpolation by PH quintics of prescribed arc length. The problem can be extended to $G^2$ data
interpolation, but the degree of the interpolating PH curve must be
elevated to $7$, since PH quintic curves do not possess enough free parameters.\\
Let us assume the complex representation and suppose we want to interpolate 
two given end points $\bfm{q}_0$, $\bfm{q}_1$, their corresponding tangent directions 
$\bfm{g}_0$, $\bfm{g}_1$, 
curvatures $\kappa_0$ and $\kappa_1$. 
Furthermore, we require that the resulting interpolant has a fixed arc length, say $L>\|\bfm{q}_1-\bfm{q}_0\|$.
Since we are looking for an interpolant $\bfm{p}$ among PH curves of degree $7$, 
there are $10$ free parameters
involved ($8$ parameters arising from the complex Bernstein coefficients \eqref{eq:preimage}
and two of them from a complex integration constant $\bfm{p}_0$).
The interpolation conditions provide 9 scalar equations. The remaining parameter could be used
for shape control or for optimization of some geometric property, but this would definitely lead to a
challenging optimization process. In order to avoid it, we will use the approach from
\cite{Farouki-quintic-arc-length-16}, i.e., we shall assume equal lengths  of the tangents of
$\bfm{p}$ at the boundary points. This reduces the number of involved
free parameters by one and gives some hope that the interpolant is fully determined already by given geometric data. The assumption is not too
restrictive and it is 
quite natural since it ensures symmetric solutions for symmetric data. \\
In \cite{Farouki-quintic-arc-length-16}, the author considered the reduction of data to the
canonical form, which has previously been used also in \cite{Farouki-Neff-Hermite-quintic-1995}.
The idea is to consider a new coordinate system which should simplify the analysis of the problem 
as much as possible. Following the above mentioned references, we can assume that the given data
is of the form $\bfm{q}_0=0$, $\bfm{q}_1=1$, $\bfm{g}_0=\exp(\compI\theta_0)$,
$\bfm{g}_1=\exp(\compI\theta_1)$ and $L>1$, where
$\theta_0,\theta_1\in(-\pi,\pi]$. The original data are transformed to the canonical one by 
an appropriate translation, rotation and scaling. Moreover, the obtained interpolant is finally pulled back to the
original coordinate system by inverse transformations.
Thus, if the above canonical data are assumed, the interpolation conditions
become
\begin{align*}
  &\frac{1}{7}
  \left(\bfm{w}_0^2
  +\bfm{w}_0\bfm{w}_1
  +\frac{3\bfm{w}_1^2+2\bfm{w}_0\bfm{w}_2}{5}
  +\frac{9\bfm{w}_1\bfm{w}_2+\bfm{w}_0\bfm{w}_3}{10}
  +\frac{3\bfm{w}_2^2+2\bfm{w}_1\bfm{w}_3}{5}
  +\bfm{w}_2\bfm{w}_3
  +\bfm{w}_3^2
  \right)-1=0,\\
  &\bfm{w}_0-d\,\exp\left(\compI\frac{1}{2}\theta_0\right)=0,\quad 
  \bfm{w}_3-d\,\exp\left(\compI\frac{1}{2}\theta_1\right)=0,\\
  &6\Im\left(\frac{\conj{\bfm{w}}_0\bfm{w}_1}{\left|\bfm{w}_0\right|^4}\right)-\kappa_0=0,\quad
  6\Im\left(\frac{\conj{\bfm{w}}_3\bfm{w}_2}{\left|\bfm{w}_3\right|^4}\right)+\kappa_1=0,\\  
  &\frac{1}{7}\left(\left|\bfm{w}_0\right|^2
  +\Re\left(\bfm{w}_0\conj{\bfm{w}}_1\right)
  +\frac{2}{5}\Re\left(\bfm{w}_0\conj{\bfm{w}}_2\right)
  +\frac{1}{10}\Re\left(\bfm{w}_0\conj{\bfm{w}}_3\right)
  +\frac{3}{5}\left|\bfm{w}_1\right|^2
  +\frac{9}{10}\Re\left(\bfm{w}_1\conj{\bfm{w}}_2\right)
  +\frac{2}{5}\Re\left(\bfm{w}_1\conj{\bfm{w}}_3\right)\right.\\
  &\left.+\frac{3}{5}\left|\bfm{w}_2\right|^2+\Re\left(\bfm{w}_0\conj{\bfm{w}}_1\right)
  +\left|\bfm{w}_3\right|^2\right)-L=0.
  \end{align*}
The first equation arises from the interpolation of two points, the
second and the third one from the interpolation of tangent directions,
the next two ensure prescribed curvatures and the last one prescribes
the arc length $L$.\\
Let us write $\bfm{w}_1=u_1+\compI v_1$, $\bfm{w}_2=u_2+\compI v_2$, 
$c_i=\cos(\theta_i/2)$, $s_i=\sin(\theta_i/2)$, $i=0,1$, replace the above equations by
their appropriate linear combinations and use some basic trigonometric identities. 
This leads to
\begin{align}
  &6u_1^2+9u_1 u_2+6u_2^2+\left(10c_0^2+10c_1^2+c_0c_1\right)d^2
  +10d\left(u_1c_0+u_2c_1\right)+4d\left(u_1c_1+u_2c_0\right)- 35(L+1)=0,\nonumber\\
  &6v_1^2+9v_1 v_2+6v_2^2+\left(10s_0^2+10s_1^2+s_0s_1\right)d^2
  +10d\left(v_1s_0+v_2s_1\right)+4d\left(v_1s_1+v_2s_0\right) - 35(L-1)=0,\nonumber\\
  &\kappa_0 d^3+6s_0u_1-6c_0v_1=0,\qquad
  \kappa_1 d^3-6s_1u_2+6c_1v_2=0,\label{eq:system}\\
&12u_1 v_1 + 9 u_2 v_1 + 9 u_1 v_2 + 12 u_2 v_2 
+\left(s_0(20c_0+c_1)+s_1(c_0+20c_1)\right)d^2\nonumber\\
&+2\left( (5 v_1 + 2 v_2) c_0 + (2v_1 + 5 v_2) c_1 + (5 u_1 + 2u_2) s_0  + 
 (2 u_1 + 5u_2) s_1\right)d=0.\nonumber
\end{align}
Since the equations arising from $G^2$ conditions are linear in $u_1$, $u_2$, $v_1$ and $v_2$, some
further reduction of the system \eqref{eq:system} is definitely possible. 
However, the analysis of the solvability for general data
seems to be extremely complicated as already guessed in \cite{Farouki-quintic-arc-length-16}.
To justify this, let us consider two particular simple examples which show that the existence of a solution
heavily depends on data.\\
Assume the data $\theta_0=\pi/2$, $\theta_1=-\pi/2$, 
$\kappa_0=\kappa_1=0$ and $L=9/8$. Plugging corresponding constants in \eqref{eq:system}
and doing some manipulations with equations reveal that $v_1=u_1$ and $v_2=-u_2$. This further
implies that either $u_2=u_1$ or $u_2=-u_1 - 5 \frac{\sqrt{2}}{6} d$.
If $u_2=u_1$, we end up with two biquadratic equations for $u_1$ and $d$, namely
\begin{equation*}
  d^2 + 8\sqrt{2}\,d\,u_1 + 18 u_1^2 - 70=0,\quad 80 d^2 + 80 \sqrt{2}\,d\,u_1 + 96 u_1^2 - 315=0.
\end{equation*}
  They actually represent a hyperbola and an ellipse shown in 
  Fig. \ref{fig:example12} (the first figure on the left) and it is easy to show
  that they do not intersect. Similarly we can check that also for $u_2=-u_1 - 5 \frac{\sqrt{2}}{6} d$
  we do not have a solution (a hyperbola and an ellipse in 
  Fig. \ref{fig:example12} (the second figure from the left)).\\
  For the second example consider the same data as  in the previous on, except that $L=2$. 
  Similar procedure as before leads to $v_1=u_1$, $v_2=-u_2$ and and again to
   $u_2=u_1$ or $u_2=-u_1 - 5 \frac{\sqrt{2}}{6} d$. If $u_2=u_1$, we end up with two biquadratic 
   equations for $u_1$ and $d$, namely
\begin{equation*}
  19 d^2 + 12\sqrt{2}\,d\,u_1 + 6 u_1^2 - 70=0,\quad 
   3 d^2 +   4\sqrt{2}\,d\,u_1 + 6 u_1^2 - 30=0,
\end{equation*}
again representing a hyperbola and an ellipse shown 
in Fig. \ref{fig:example12} (the third figure from the left).  It is clearly seen that they
intersect in four points, thus the system of nonlinear equations \eqref{eq:system}
has four solutions. The case $u_2=-u_1 - 5 \frac{\sqrt{2}}{6} d$ again implies a hyperbola and an ellipse 
with no intersections (Fig. \ref{fig:example12}, the last figure in the row). It is clear that changing also the angles
$\theta_i$ and curvatures $\kappa_i$, $i=0,1$, would imply even more complicate examples with solutions relying heavily on the data.
  
\section{Interpolation of circular arc data}\label{sec:circular_arc}

The analysis of the system of nonlinear equations
determining PH curve of degree 7, which interpolates given $G^2$ data with a prescribed arc length,
is in general highly nontrivial tack as can be seen from the examples in the previous section. 
In this section we shall consider the same interpolation problem, but for circular arc data, 
i.e., the data arising from a circular arc. Even if this simplification is considered, there seems to be no results
available in the literature, so the problem is worth studying.
A similar problem was considered in \cite{Farouki-et-al-2021-PH7-clothoid}
where the authors considered the arc length preserving approximations 
of monotone clothoid segments.\\
Let the data be sampled from the circular arc with its inner angle
equal to $2\alpha$. The canonical position and some elementary 
geometry imply $\alpha=\theta_0=-\theta_1$,
$\kappa_0=\kappa_1=-2\sin\alpha$, 
$L=\alpha\csc\alpha$,
the radius of the arc equals to $1/(2\sin\alpha)$ and its centre is
$(1/2,-1/2\cot\alpha)^T$.
Note that in the following we will use two additional standard 
trigonometric functions $\csc=1/\sin$ and $\sec=1/\cos$.
We shall further assume that
$0<\alpha\leq \pi/2$, since for practical applications it is enough to construct good approximations
of circular arcs up to the semicircle. It also seems that similar analysis as in the following could be done
for $\pi/2<\alpha<\pi$, too.\\
The circular data are first used to determine constants in the nonlinear system \eqref{eq:system}.
The third and the fourth equation are then solved on 
$v_1$ and $v_2$, i.e., 
\begin{equation}\label{eq:v1v2}
  v_1=\frac{1}{3} \tan\left(\frac{\alpha}{2}\right)
  \left(3 u_1 - 2 d^3 \cos\left(\frac{\alpha}{2}\right)\right),\quad
  v_2=-\frac{1}{3} \tan\left(\frac{\alpha}{2}\right)
  \left(3 u_2  - 2 d^3 \cos\left(\frac{\alpha}{2}\right)\right).
\end{equation}
Combining this with the fifth equation leads to
\begin{equation*}
  (u_1 - u_2) \left(6 (u_1 + u_2) - d (d^2-10) \cos\left(\frac{\alpha}{2}\right)\right)=0. 
  \end{equation*}
We obviously have two possibilities,
$u_2=-u_1+\frac{1}{6}d(d^2-10)\cos\left(\frac{\alpha}{2}\right)$ or
$u_2=u_1$. Note that \eqref{eq:system} implies that 
if $(d,u_1,v_1,u_2,v_2)$ is a solution, 
then also $(-d,-u_1,-v_1,-u_2,-v_2)$ is a solution 
which, by \eqref{eq:p_coeff}, provides the same interpolant $\bfm{p}$,
so it is enough to consider solutions with $d>0$ only.

\subsection{The case
$u_2=-u_1+\frac{1}{6}d(d^2-10)\cos\left(\frac{\alpha}{2}\right)$}
\label{sub:nonsymetric}

We will prove that there are no real solutions in this case. 
Some particular linear combinations of the first two equations in
\eqref{eq:system} lead to two nonlinear equations
\begin{align}
&18(4\cos\alpha-3)u_1^2
+3d(d^2-10)\left(\cos\left(\frac{\alpha}{2}\right)
-2\cos\left(\frac{3\alpha}{2}\right)\right)u_1\nonumber\\
&+\cos^2\left(\frac{\alpha}{2}\right)
\left(-3d^6+18 d^4-34d^2-420
+4 d^2 (d^4-6d^2+30)\cos\alpha\right)=0,\label{eq:sys_nonsym}\\
&36\left(3\sin\left(\frac{3\alpha}{2}\right)
-11\sin\left(\frac{\alpha}{2}\right)\right)u_1^2
+6\sin\left(\frac{\alpha}{2}\right)
\left(5\cos\left(\frac{\alpha}{2}\right)
-3\cos\left(\frac{3\alpha}{2}\right)\right)
d (d^2-10)u_1\nonumber\\
&+\left(840\alpha - 8 d^2 (d^4-6d^2+30)\sin\alpha 
+ d^2 ( 3d^4-18d^2+34)\sin\left(2\alpha\right)\right)
   \cos\left(\frac{\alpha}{2}\right)=0.\nonumber
\end{align}
Fortunately, the resultant of the polynomials on the left hand side
of \eqref{eq:sys_nonsym} with respect to $u_1$
simplifies to
\begin{equation}\label{eq:r_nonsym}
    r(d)=\left(16\sin^3\alpha\,d^2
    +30((3\cos\alpha-4)\sin\alpha+\alpha(4\cos\alpha-3))\right)^2
\end{equation}
and the candidates for solutions of \eqref{eq:sys_nonsym} with 
positive $d$ are positive zeros of $r$.

\begin{lem}\label{lem:resultant_nonsym}
  For $\alpha\in(0,\pi/2]$, function $r$ has precisely one (double) positive zero
  $$
    d_1=\frac{\sqrt{30}}{4}\sqrt{\frac{\alpha(3-4\cos\alpha)
    +(4-3\cos\alpha)\sin\alpha}{\sin^3\alpha}}>\frac{5}{2}.
  $$
\end{lem}
\begin{pf}
  By \eqref{eq:r_nonsym} (double) zeros of $r$ are $\pm d_1$.
  The result of the lemma will follow if we prove that $f(\alpha)>0$
  on $(0,\pi/2]$, where
  $$
    f(\alpha)=\alpha(3-4\cos\alpha)
    +(4-3\cos\alpha)\sin\alpha-\frac{10}{3}\sin^3\alpha.
  $$
  Quite clearly $f(0)=0$ and 
  $
  f'(\alpha)=\sin\alpha(4\alpha+6\sin\alpha
  -5\sin2\alpha)\geq 6\sin^2\alpha\,(1-\cos\alpha)>0,
  $
  where we have used the known fact that $\alpha>\sin\alpha$ for
  $\alpha>0$. Consequently $f$ is positive on $(0,\pi/2]$
  and the proof is completed.
\qed\end{pf}
Note that the second equation in \eqref{eq:sys_nonsym} is quadratic
in $u_1$ with the discriminant 
\begin{equation}\label{eq:discriminant_nonsym}
  504\left(2\cos\left(\frac{\alpha}{2}\right)
  +6\sin\left(\frac{\alpha}{2}\right)\sin\alpha\right)
  \sin\left(\frac{\alpha}{2}\right)f_\alpha(d),
\end{equation}
where
$$
f_\alpha(d)=960 \alpha + 
 -8d^2(d^4-4d^2+20)\sin\alpha
 +d^2(3d^4-12d^2-4) \sin2\alpha.
$$
\begin{lem}\label{lem:f_alpha_nonsym}
  Function $f_\alpha$ is negative on $[5/2,\infty)$ for all
  $\alpha\in(0,\pi/2]$.
\end{lem}
\begin{pf}
  The idea of the proof is similar as in the proof of
  Lemma \ref{lem:resultant_nonsym}. Let 
  $$
  g(\alpha)=f_\alpha(5/2)=
  960 \alpha -\frac{13625 \sin\alpha}{8}+\frac{15275}{64} \sin2\alpha.
  $$
  Obviously $g(0)=0$ and 
  $g'(\alpha)=\tfrac{5}{32}(6110\cos^2\alpha-10900\cos\alpha+3089)$.
  By solving a simple quadratic equation one can conclude that
  $\alpha_0\approx 1.2069$ is the unique zero of $g'$ on $(0,\pi/2]$
  implying the local minimum $g(\alpha_0)\approx-274.2089$.
  Since also $g(\pi/2)\approx-195.1605$, function $g$ must be 
  negative on $(0,\pi/2]$. Furthermore, zeros of $f_\alpha'$
  are $d_1=0$, 
  $$d_{2,3}=
  \pm \sqrt{\frac{16-12\cos\alpha-h(\alpha)}{12-9\cos\alpha}},
  \quad 
  d_{4,5}=
  \pm \sqrt{\frac{16-12\cos\alpha+h(\alpha)}{12-9\cos\alpha}},
  $$
  where
  $$
    h(\alpha)=\sqrt{-614+288 \cos\alpha+90 \cos2\alpha}.
  $$
 Obviously $h^2<0$, and $d_{2,3,4,5}$ are complex thus $f_\alpha$
 must be monotone on $[5/2,\infty)$. Since
 $f_\alpha''(0)=-16\sin\alpha(20\sin\alpha+\cos\alpha)<0$,
 $f_\alpha$ is decreasing which together with $f_\alpha(5/2)<0$
 implies the result of the lemma.
\qed
\end{pf}
From Lemma \ref{lem:resultant_nonsym},
Lemma \ref{lem:f_alpha_nonsym} and \eqref{eq:discriminant_nonsym}
now follows that the system of nonlinear
equations \eqref{eq:sys_nonsym} has no real solutions for any
$\alpha\in(0,\pi/2]$.

\subsection{The case $u_1=u_2$}\label{sub:symmetric}
\begin{figure}[htb]
\centering
\includegraphics[width=0.5\textwidth]{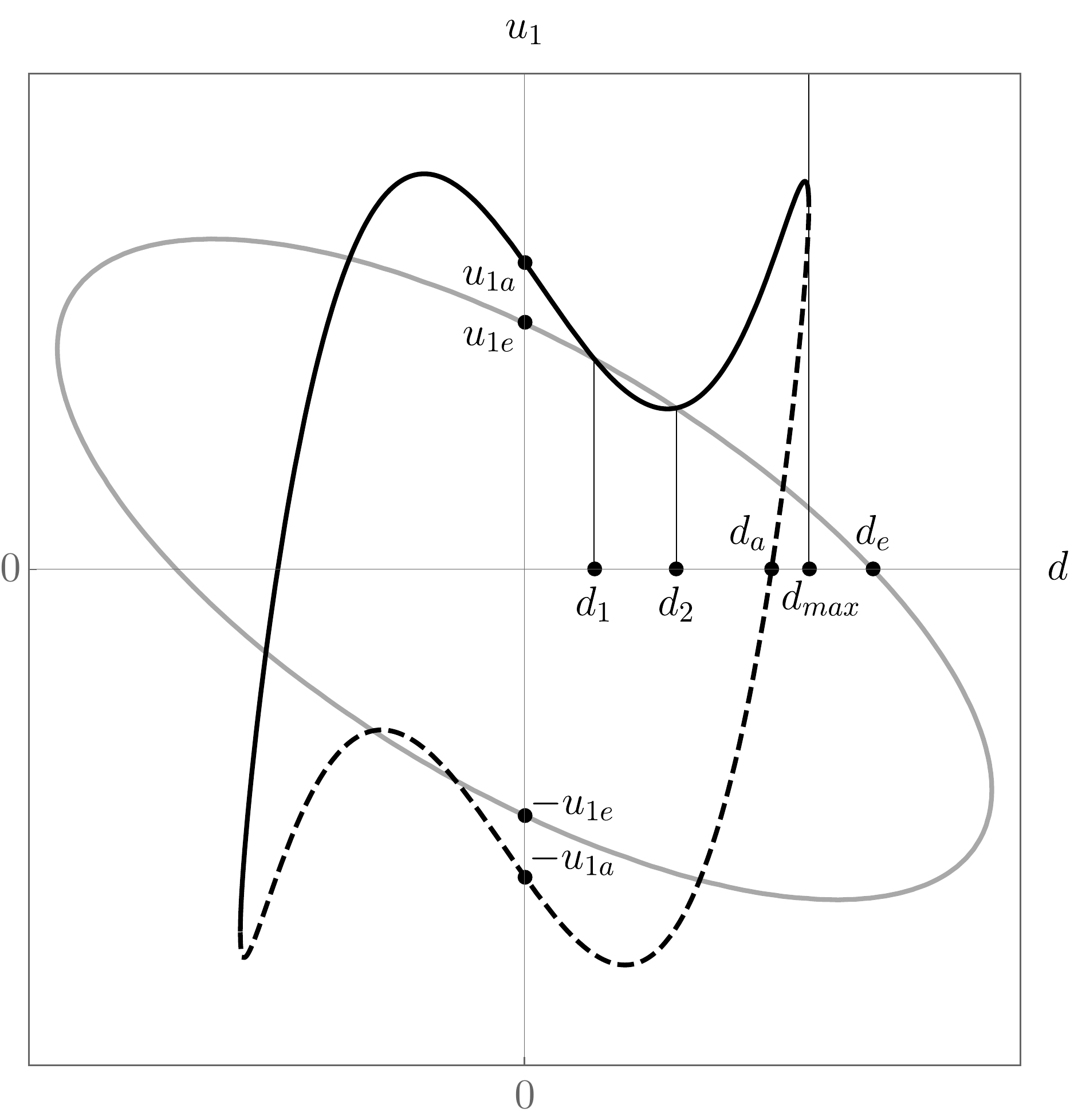}
\caption{Curves given by equations \eqref{eq:curve1} (gray, an ellipse) and \eqref{eq:curve2}
(black solid and black doted, an algebraic curve of degree $6$) 
with quantities used in the proof of Lemma \ref{lem:at_least_two_solutions}.}
\label{fig:curves12}
\end{figure}

Let us consider now the symmetric case, i.e., $u_1=u_2$.
Equations \eqref{eq:v1v2} imply $v_2=-v_1$ and some particular linear combinations of the first two equations in 
\eqref{eq:system} lead to
 \begin{align}
  &3(1+\cos\alpha) d^2 + 8\cos\left(\frac{\alpha}{2}\right)d u_1 + 6u_1^2 - 
  10(1+\alpha\csc\alpha)=0,\label{eq:curve1}\\
  &4 d^6
  - 24 d^4 
 +57 d^2
  -12\sec\left(\frac{\alpha}{2}\right)  d (d^2-3) u_1
  +9\sec^2\left(\frac{\alpha}{2}\right)u_1^2
  +105\csc^2\left(\frac{\alpha}{2}\right)\left(1-\alpha\csc\alpha\right)=0,\label{eq:curve2}
\end{align}
the system of two nonlinear equations for $d$ and $u_1$.
The first equation again represents an ellipse, while
the second one is much more complicated algebraic
curve of degree $6$ (see Fig. \ref{fig:curves12}). But it is
quadratic in $u_1$ and we can detect its two branches 
\begin{equation}\label{eq:u1(d)}
    u_1(d)=\frac{1}{3}\cos\left(\frac{\alpha}{2}\right)
    \left(2d^3-6d\pm
    \sqrt{21}\sqrt{5\csc^2\left(\frac{\alpha}{2}\right)
    \left(\alpha\csc{\alpha}-1\right)-d^2}\right)
  \end{equation}
(the plus sign is the black solid part of the curve and the 
minus sign is the black dashed part in Fig. \ref{fig:curves12}). 
It is also easy to see that
it is a closed curve with 
$d\in[-d_{max},d_{max}]$
where
\begin{equation}\label{eq:dmax}
    d_{max}=\csc\left(\frac{\alpha}{2}\right)\sqrt{5\left(\alpha\csc\alpha-1\right)}.
  \end{equation}
\begin{lem}\label{lem:dmax_bound}
  For $\alpha\in(0,\pi/2]$ the inequality 
  $d_{max}>\sqrt{10/3}$ holds true.
\end{lem}
\begin{pf}
  It is easy to see that the inequality from the lemma 
  is equivalent to $f(\alpha)>0$, where 
  $f(\alpha)=3(\alpha-\sin\alpha)-\sin\alpha(1-\cos\alpha)$.
  Since $f(0)=0$ and 
  $f'(\alpha)=8\sin^4\left(\frac{\alpha}{2}\right)>0$, the proof of the
  lemma is complete.
  \qed
\end{pf}
The following lemma guarantees the existence of the solution of the considered interpolation problem.
\begin{lem}\label{lem:at_least_two_solutions}
  The system of nonlinear equations \eqref{eq:curve1}, \eqref{eq:curve2} has at least two
  real solutions with $d>0$ for any $\alpha\in(0,\pi/2]$.
\end{lem}
\begin{pf}
  During the proof we will refer to Fig. \ref{fig:curves12}. 
  Consider the $(d,u_1)$ plane. 
  Let $u_{1e}$ and $u_{1a}$ be the positive intersections of \eqref{eq:curve1} and \eqref{eq:curve2}
  with $d=0$, respectively, and similarly, let $d_e$ and $d_a$ be the
  positive intersections of the same
  curves with $u_1=0$, respectively. If we show that
  \begin{equation}\label{neq:u11u22d1d2}
    (u_{1a}-u_{1e})(d_a-d_e)<0,
  \end{equation}
  then curves must intersect in the first quadrant of the chosen coordinate system. 
  But due to the symmetric properties they must also intersect in the fourth quadrant
  and the result of the lemma will follow.\\
  Quite clearly, for \eqref{neq:u11u22d1d2} it is enough to see $u_{1a}>u_{1e}$ and $d_e>d_a$.
  Let us start by proving the first inequality. If $d=0$, then \eqref{eq:curve1} and \eqref{eq:u1(d)} imply
  \begin{align*}
    u_{1e}&=\sqrt{\frac{5}{3}}\sqrt{\alpha\csc\alpha+1},\quad
    u_{1a}=\sqrt{\frac{35}{3}}\cot\left(\frac{\alpha}{2}\right)\sqrt{\alpha\csc\alpha-1}.    
  \end{align*}
  Some straightforward calculations reveal that $u_{1a}>u_{1e}$ is equivalent to
  $f(\alpha):=\cos\alpha\left(4\alpha - 3\sin\alpha\right)- 4\sin\alpha+3\alpha>0$ on
  $(0,\pi/2]$.  Since $f'(\alpha)=2\sin\alpha(3\sin\alpha-2\alpha)$, $f$ is strictly increasing on
  $(0,\alpha^*)$ and strictly decreasing on $(\alpha^*,\pi/2]$, where $\alpha^*\in(\pi/4,\pi/2)$.
  But $f(0)=0$, $f(\pi/2)=(3\pi-8)/2>0$, and the conclusion $f>0$ on $(0,\pi/2]$ follows.\\
  For the second inequality observe that $d_{max}>d_a$ and it is enough to see that $d_e\geq d_{max}$.
  Inserting $u_1=0$ in \eqref{eq:curve1} and considering 
  \eqref{eq:dmax} lead us to show that
  $g(\alpha):=\cos\alpha\sin\alpha+2\sin\alpha-2\alpha\cos\alpha-\alpha\geq 0$ on $(0,\pi/2]$.
  But this follows immediately from $g(0)=0$ and $g'(\alpha)=2\sin\alpha(\alpha-\sin\alpha)>0$
  on $(0,\pi/2]$. 
\qed
\end{pf}

\begin{figure}[htb]
\centering
\includegraphics[width=1\textwidth]{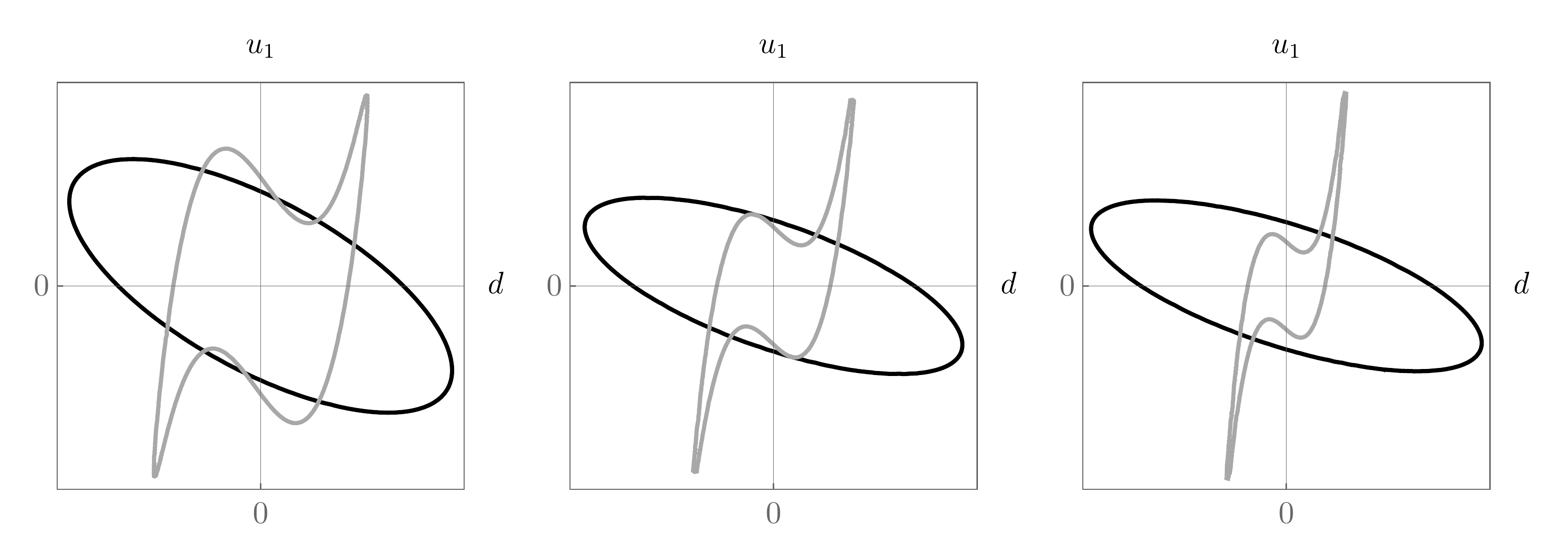}
\caption{Four solutions with positive $d$ for $\alpha<\alpha_{crit}\approx
2.2337$ (left) which degenerate to three of them
for $\alpha=\alpha_{crit}$ (middle) and transform to two solutions
for $\alpha>\alpha_{crit}$ (right).}
\label{fig:general_cases}
\end{figure}

\begin{rem}
  Identical proof can be done for the case $\alpha\in (0,\alpha_{max})$, where 
  $\alpha_{max}\approx 2.0682$ is the first positive zero of $f$ defined in the proof of the previous lemma.
\end{rem}
From the previous lemma it actually follows that the
considered system of nonlinear equations has an even
number of solutions with positive $d$ for any
$\alpha\in(0,\pi/2]$. Numerical
examples reveal that there might be four of them as
indicated in Fig. \ref{fig:curves12}. However,
it seems quite difficult to prove the 
existence of other two solutions ($d_1$ and $d_2$ on 
Fig. \ref{fig:curves12}) in general, since they tend to each other for small
$\alpha$ and they disappear for 
$\alpha>\alpha_{crit}\approx 2.2337$, i.e., the critical value of $\alpha$ which is
a solution of \eqref{eq:curve1}, \eqref{eq:curve2} their zero
Jacobian (see Fig. \ref{fig:general_cases}).
In the following, we will provide an easy (numerical)
procedure
to check the existence of four solutions with $d>0$ for
a particular $\alpha\in(0,\pi/2]$ and prove their existence for
$\alpha$ small enough. The PH interpolant of degree 
seven arising from the positive solution $d_j$, $j=1,2,3,4$, where
$d_1<d_2<d_3<d_4$, will
be denoted by $\bfm{p}_j$.\\
Let us first transform the system of nonlinear equations 
\eqref{eq:curve1} and \eqref{eq:curve2} to a more appropriate one
for the analysis. This can be done by using a Gr\"oebner basis
(\cite{Adams-Loustaunau-Groebner-1994}) with respect to a particular ordering of the unknowns.
The system of nonlinear equations then reads as

\begin{align}
  p_1(d)
  &=-32 \sin^6\alpha\, d^{12}
  +256 \sin ^6\alpha\,d^{10}
  -1184 \sin^6\alpha\,d^{8}\nonumber\\
  &-96 \sin^3\alpha(-40 \alpha +9 \sin\alpha
  +20 \sin2\alpha+7\sin3\alpha
  -30 \alpha \cos\alpha)\,d^{6}\nonumber\\
  &+96 \sin^3\alpha(-160\alpha +99 \sin\alpha
  +80 \sin2\alpha+7 \sin3\alpha
  -120\alpha \cos\alpha)\,d^{4}\nonumber\\
  &+13440 (\alpha -\sin\alpha)\sin^5\alpha
  \csc^2\left(\frac{\alpha }{2}\right)\,d^{2}\nonumber\\
  &-1800 (6\alpha +8\alpha\cos\alpha
  -2\sin\alpha(3\cos\alpha+4))^2=0,\label{eq:p1}\\
  p_2(u_1,d)&=
  24d\left(d^2-2\right)\sec\left(\frac{\alpha }{2}\right)u_1
  +210 \csc^2\left(\frac{\alpha }{2}\right) 
  (\alpha  \csc\alpha-1)\nonumber\\
  &+\sec^2\left(\frac{\alpha }{2}\right) 
  \left(-30 \alpha  \csc\alpha+9 d^2 \cos\alpha
  +9d^2-30\right)
  -8 d^6+48 d^4-114d^2=0.\label{eq:p2}
  \end{align}
Let us analyze the polynomial $p_1$ first. Since it is even of degree
$12$, we can reduce its degree to $6$ by introducing
$p(x)=p_1(\sqrt{x})$. 
A computer algebra system reveals that 
\begin{equation}\label{eq:p_expansion}
  p(x)=-32q(x)\alpha^6+{\mathcal O}(\alpha^8),\quad 
  q(x)=x^6-8 x^5+37 x^4-134 x^3+284 x^2-280 x+100,
\end{equation}
and real zeros of $q$ are
\begin{equation}\label{eq:xi}
x_1=x_2=1\,\quad x_3\approx 2.1842,\quad x_4\approx 3.2872.
\end{equation}
Note that $x_3$ and $x_4$ can also be written in radicals 
since they are zeros of some quartic polynomial but expressions are to complicated to be
given here explicitly.

\begin{lem}\label{lem:four_solutions_p}
  For any $\alpha\in(0,\pi/2]$, 
  the polynomial $p$ has at most four positive 
  zeros on $(0,d_{max}^2)$. If $p(y_i)$, $i=0,1,\dots,4$,
  where $y_0=0$, $y_1=1$, $y_2=x_3$, $y_3=x_4$ and
  $y_4=d_{max}^2$, are of alternating signs,
  then $p$ has precisely four positive zeros.
\end{lem}
\begin{pf}
  First observe that the fourth derivative of $p$
  is of particularly simple form, namely
  $p^{(iv)}(x)=-768 \sin^6\alpha\left(15 x^2-40 x+37\right)$. Quite clearly it is negative and 
  consequently $p$ has at most four real zeros.
  By Lemma \ref{lem:dmax_bound} and \eqref{eq:xi} we have
  $y_0<y_1<y_2<y_3<y_4$ and $p(y_i)p(y_{i+1})<0$, $i=0,1,2,3$, 
  implies precisely four zeros due to the continuity of $p$.
  \qed
\end{pf}
We are now ready to prove the main theorem of this paper.
\begin{thm}
  If $\alpha\in(0,\pi/2]$ then the number of real solutions of the nonlinear system
  \eqref{eq:curve1}, \eqref{eq:curve2} with $d>0$ is
  the same as the number of positive zeros of $p$.
\end{thm}
\begin{pf}
  Since nonlinear system \eqref{eq:curve1}, \eqref{eq:curve2}
  is equivalent to the system \eqref{eq:p1}, \eqref{eq:p2},
  the only candidates for real solutions with $d>0$ are,
  by Lemma \ref{lem:four_solutions_p}, positive zeros of $p$.
  Thus we have to prove that each positive zero of $p$ implies 
  the unique real solution of the system \eqref{eq:p1}, \eqref{eq:p2}.
  Let $z$ be a positive zero of $p$. Then $d_z=\sqrt{z}$
  is a positive zero of $p_1$ and the solution of 
  $p_2(u_1,d_z)=0$ on $u_1$ provides the desired solution of
  the system of nonlinear equations. But $p_2(u_1,\cdot)$ is
  a linear polynomial and it remains to prove that its leading
  coefficient does not vanish at $d_z$. It is equivalent to
  verifying that $z\neq 0,2$, or equivalently 
  \begin{equation*}
      p(0)=-7200f_1(\alpha)^2<0,\quad
      p(2)=-32f_2(\alpha)^2<0,
  \end{equation*}
  where
  \begin{align*}
    f_1(\alpha)&=\cos\alpha\left(4\alpha - 3\sin\alpha\right)- 4\sin\alpha+3\alpha,\\
    f_2(\alpha)&=45 \alpha -2\sin ^3\alpha
    -66 \sin\alpha+60\alpha\cos\alpha
    +6 \sin\alpha\cos^2\alpha-45\sin\alpha \cos\alpha.
    \end{align*}
  Thus it is enough to show that $f_1,f_2>0$ on $(0,\pi/2]$. 
  The inequality $f_1>0$ follows from the fact that $f_1=f$ from the
  proof of Lemma \ref{lem:at_least_two_solutions}. To confirm $f_2>0$,
  observe that $f_2(0)=0$ and
  $f_2'(\alpha)=-6\sin\alpha f_3(\alpha)$, where $f_3(\alpha)=10\alpha-15\sin\alpha+4\sin\alpha\cos\alpha$.
  Since $f_3'(\alpha)=6-15\cos\alpha+8\cos^2\alpha$,
  $f_3'$ has precisely one zero on $[0,\pi/2]$ and consequently
  $f_3$ has at most two zeros there. Since $f_3(0)=0$,
  $f_3(\alpha)=-\alpha+{\mathcal O}(\alpha^3)$ and 
  $f_3(\pi/2)=5(\pi-3)>0$, $f_3$ has the unique zero
  $\alpha_0\in(0,\pi/2]$.
  Thus $f_2$ is increasing on $(0,\alpha_0)$ and decreasing 
  on $(\alpha_0,\pi/2]$. Since $f_2(0)=0$ and $f_2(\pi/2)=(45\pi-136)/2>0$, function $f_2$ must be positive on
  $(0,\pi/2]$ and the result of the theorem follows.
  \qed 
\end{pf}
Previous theorem provides an efficient and easy way to check the number
polynomial parametric approximants interpolating $G^2$ data and
an arc length arising
from a circular arc given by an inner angle $2\alpha$. 
For some practically important angles $\alpha$, such as 
$\alpha=\pi/2,\pi/3,\pi/4,\pi/8,\dots$, the direct application of 
Lemma \ref{lem:four_solutions_p} confirmed the existence 
of four zeros of $p$, except for $\alpha=\pi/2$, where we had
to replace $y_2=x_3$ by $y_2=2$. A direct formal proof that precisely
four solutions exist for any $\alpha\in(0,\pi/2]$ seems to be
quite a difficult task, since the analysis of symbolic expressions
involving combinations of algebraic and trigonometric terms 
in Lemma \ref{lem:four_solutions_p} would
be needed. However, if $\alpha$ is small enough, the expansion
\eqref{eq:p_expansion} enables us to prove the existence of four
solutions in general. This will be done in the following section.

\section{Asymptotic analysis}\label{sec:asymptotic_analysis}

Let us now consider $\alpha$ small enough. Using \eqref{eq:p_expansion}
and considering some additional terms in the expansion, we get
\begin{align*}
  p(y_0)&=-3200\alpha^6+{\mathcal O}(\alpha^8),\quad 
  p(y_1)=64\alpha^{10}+{\mathcal O}(\alpha^{12}),\quad
  p(y_2)\approx-52.3867\alpha^{8}+{\mathcal O}(\alpha^{10}),\\
  p(y_3)&\approx 2292.89\alpha^{8}+{\mathcal O}(\alpha^{10}),\quad
  p(y_4)=-\frac{156800}{729}\alpha^{6}+{\mathcal O}(\alpha^{8}).
\end{align*}
Consequently, $p$ has four positive zeros by Lemma \ref{lem:four_solutions_p}. The leading terms 
constants can be written also in a closed form, 
thus their numerical values can be computed with arbitrary precision.\\
In the following we will find asymptotic expansions of positive
zeros of $p$ which provide asymptotic expansions of real
solutions of the system of nonlinear equations \eqref{eq:curve1},
\eqref{eq:curve2}.
Let $z_i$, $i=1,2,3,4$, be a positive zero of $p$. Then 
\eqref{eq:p_expansion} suggests the expansion of $z_i$ as
\begin{equation*}
  z_i=x_{i}+\sum_{j=1}^\infty c_{i,j}\alpha^j,\quad i=1,2,3,4.
\end{equation*}
Constants $c_{i,j}$ can now be found as a solution of the system of
equations for $c_{i,j}$ arising from the condition that terms in the expansion of $p(z_i)$ vanish for all $\alpha$. Let us demonstrate
the procedure for $i=2$, since the solution $z_2$ will later 
turn out as the most appropriate one. The expansion of $p(z_2)$
reads as
\begin{align*}
    p(z_2)&=-1248c_{2,1}^2\alpha^8 
    +64c_{2,1} \left(23 c_{2,1}^2-39 c_{2,2}-2\right)\alpha^9\\
    &-32\left(12 c_{2,1}^4-3 \left(46 c_{2,2}+9\right) 
    c_{2,1}^2+78 c_{2,3} c_{2,1}+39c_{2,2}^2+4 c_{2,2}-2\right)
     \alpha ^{10}+\dots
\end{align*}
The requirement that the coefficients at $\alpha^j$, $j=8,9,10$,
vanish, leads to the triangular system of nonlinear equations
with solutions 
$c_{2,1}=0,\quad c_{2,2}=(-2\pm\sqrt{82})/39$. Since $z_2>1$, we must
take $c_2=(-2+\sqrt{82})/39$. Considering more terms in the expansion,
we can similarly compute additional constants $c_{2,j}$ 
but we will skip the details. Recall that $d_2=\sqrt{z_2}$, so 
the asymptotic expansion of $d_2$ is
\begin{equation}\label{eq:d2asym}
  d_2=1-\frac{1}{78} \left(2-\sqrt{82}\right)\alpha^2
  +\frac{\left(37966+10579 \sqrt{82}\right)}{19456632}\alpha^4
  +\cdots
\end{equation}
Together with \eqref{eq:p2} we get the asymptotic expansions
\begin{equation}\label{eq:u12asym}
  u_{1,2}=u_{2,2}=1
  +\frac{1}{312} \left(73-4 \sqrt{82}\right)\alpha^2
  -\frac{\left(3071515-636632\sqrt{82}\right)}{311306112}\alpha^4
  +\cdots
\end{equation}
and finally from \eqref{eq:v1v2} also
\begin{equation}\label{eq:v1v2asym}
  v_{1,2}=-v_{2,2}=\frac{\alpha }{6}
  +\frac{\left(371-36\sqrt{82}\right)}{1872}\alpha^3
  -\frac{\left(8660963-1179080\sqrt{82}\right)}{1037687040}\alpha^5
  +\cdots
\end{equation}
Similarly we compute asymptotic expansions for other three solutions
$d_1$, $d_3$, $d_4$ and consequently also expansions for
$u_{1,j}$, $u_{2,j}$, $v_{1,j}$ and $v_{2,j}$, $j=1,3,4$. 
Either in non-asymptotic or in asymptotic approach we obtain 
several solutions. In the next section we will provide some
suggestions how to choose the most appropriate one.

\section{Solution selection}\label{sec:solution_selection}

Multiple solutions are regularly observed fact when one is dealing with interpolation by PH curves. 
Usually, some of them are more appropriate for applications
(without undesirable loops, e.g.) than the other ones.
This was observed already in the early papers dealing
with interpolation by PH curves 
(\cite{Albrecht-Farouki-C2-PH-1996}, \cite{Farouki-Neff-Hermite-quintic-1995}). 
There are several suggestions how to choose the most appropriate solution,
but non of them can be considered as a universal one.
Quite standard measures of fairness is
the absolute rotation index (\cite[p. 532]{Farouki-PH-book-2008}),
which is defined as
\begin{equation}\label{eq:rotation_index}
  R_{abs}=\int_{0}^1 \left|\kappa(t)\right|\|\bfm{p}'(t)\|{\rm d}t.
\end{equation}
It was successfully used in \cite{Farouki-quintic-arc-length-16} to identify more appropriate solutions.
We can use the same criterion here for the general interpolation of $G^2$ data. However, for the circular arc data it seems reasonable to observe the deviation of the curvature of the interpolant from the
(constant) curvature of the corresponding circular arc in $L^2$ norm. Since the curvature of the
circular arc in the chosen canonical position is $-2\sin\alpha$, 
the error becomes
\begin{equation}\label{eq:curvature_error}
  E_\kappa=\int_{0}^1 \left(\kappa(t)+2\sin\alpha\right)^2{\rm d}t.
\end{equation}
Note that the (numerical) evaluation of $E_\kappa$ for a PH curve is particularly simple since its curvature $\kappa$ is 
a rational function.\\
For the asymptotic case explained in the previous section, it is
promising to choose the solution which provides a curve with
the best approximation properties, such as the minimal 
Hausdorff distance. Since the approximation of a circular arc is
considered, one can use the radial distance $d_{rad}$
as the error measure (\cite{Degen-92-Best-approximation}). 
It is a special type of parametric distance considered
in \cite{Lyche-Morken-94-Metric} and later in
\cite{Jaklic-Kozak-2018-best-circle} where the authors 
have proved that it actually coincides with the Hausdorff distance
in the case of circular arc approximation. 
Let $\bfm{p}=(\bfm{p}_x,\bfm{p}_y)^T$ 
be a PH curve of degree seven approximating
the circular arc $\bfm{c}$ given by some small inner angle $2\alpha$
in the canonical position.
Then the radial distance is defined as
\begin{equation*}
  d_{rad}(\bfm{p};\alpha)=\max_{t\in[0,1]}
  \left|\sqrt{\left(\bfm{p}_x(t)-\frac{1}{2}\right)^2+
  \left(\bfm{p}_y(t)+\frac{1}{2}\cot\alpha
  \right)^2}
  -\frac{1}{2\sin\alpha}\right|.
  \end{equation*}
i.e., the distance between the point $\bfm{p}(t)$ and the intersection
of the line passing through the centre of the circular arc
$(1/2,-1/2\cot\alpha)^T$ and $\bfm{p}(t)$ with the circular arc.
Using the asymptotic expansions 
\eqref{eq:d2asym}--\eqref{eq:v1v2asym} we can derive
\begin{equation}\label{eq:radial_distance}
    d_{rad}(\bfm{p};\alpha)=
    c\,\alpha^r+{\mathcal O}\left(\alpha^{r+1}\right),
\end{equation}
where $c$ is some positive constant and $r\in\NN$ is the
asymptotic approximation order. Since $\bfm{p}$ interpolates
two points, two tangent directions, two curvatures and an arc length,
the expected approximation order is $7$. Indeed, for the solution 
$d=d_2$ which implies the PH interpolant $\bfm{p}_2$ we have
\begin{equation}\label{eq:order2}
  d_{rad}(\bfm{p}_2;\alpha)=
  \frac{47773-5264 \sqrt{82}}{318898944}\alpha^7
  +{\mathcal O}\left(\alpha^{8}\right)\approx
  3.3068\times10^{-7}\,\alpha^7
  +{\mathcal O}\left(\alpha^{9}\right),
\end{equation}
while $d_1$, $d_3$ and $d_4$ imply interpolants 
$\bfm{p}_1$, $\bfm{p}_3$ and $\bfm{p}_4$, respectively, 
with inferior leading term constant or much lower approximation order.
More precisely, 
\begin{align}
  d_{rad}(\bfm{p}_1;\alpha)&=
  \frac{47773+5264 \sqrt{82}}{318898944}\alpha^7
  +{\mathcal O}\left(\alpha^{8}\right)\approx
  2.9928\times10^{-4}\,\alpha^7
  +{\mathcal O}\left(\alpha^{9}\right),\label{eq:order1}\\
  d_{rad}(\bfm{p}_3;\alpha)&=0.0173\alpha +{\mathcal O}\left(\alpha^{3}\right),\quad 
  d_{rad}(\bfm{p}_4;\alpha)=0.1246\alpha +{\mathcal O}\left(\alpha^{3}\right).\label{eq:order34}
\end{align}

\section{Numerical examples}\label{sec:examples}

\begin{figure}[htb]
\centering
\includegraphics[width=0.75\textwidth]{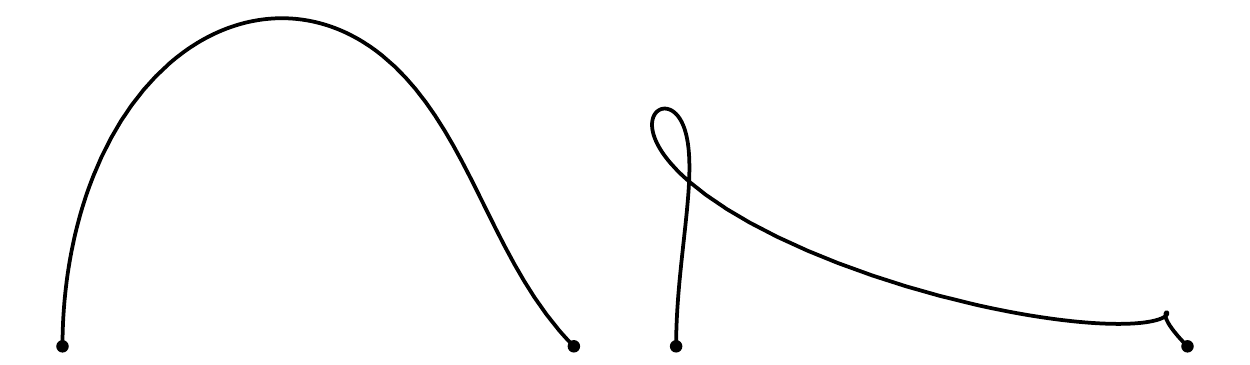}
\caption{Two solution curves interpolating data \eqref{ex:1}.}
\label{fig:example1}
\end{figure}

\begin{figure}[htb]
\centering
\includegraphics[width=0.8\textwidth]{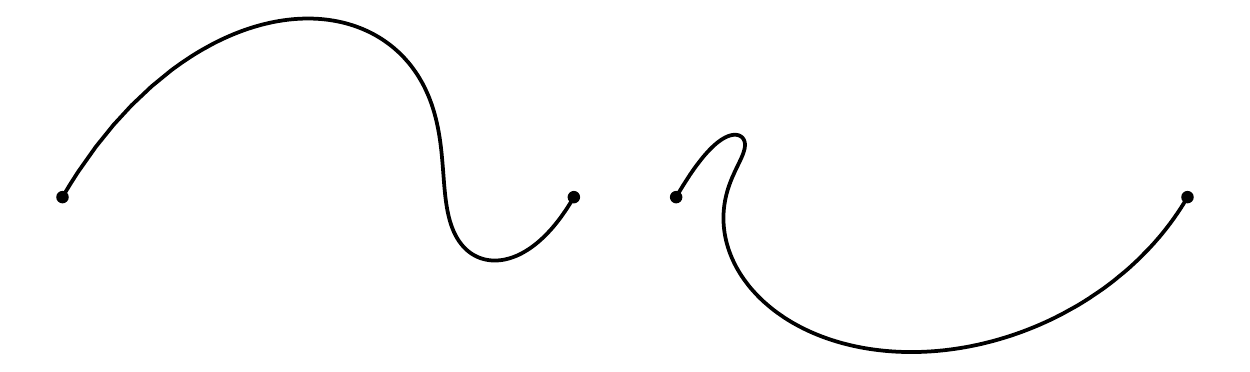}
\caption{Plots of interpolants of data \eqref{ex:2}.}
\label{fig:example2}
\end{figure}

\begin{figure}[htb]
  \centering
  \includegraphics[width=0.8\textwidth]{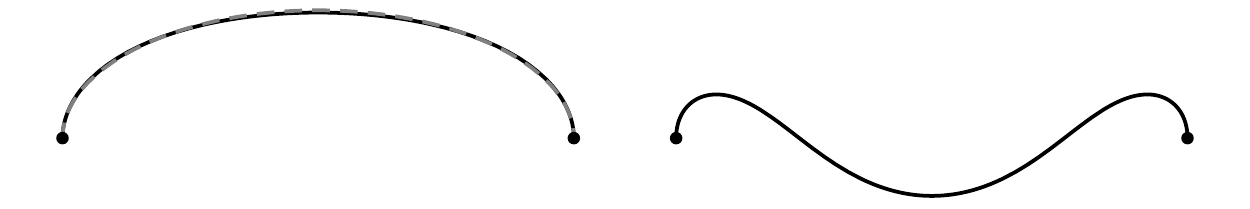}
  \caption{Two formally admissible interpolants for data \eqref{ex:3}. Almost invisible dashed gray
  curve on the left is the elliptical arc.}
  \label{fig:example3}
  \end{figure}

\begin{figure}[!htb]
\begin{tabular}{c}
 \minipage{0.25\textwidth}
    \centering
    \includegraphics[width=0.8\linewidth]{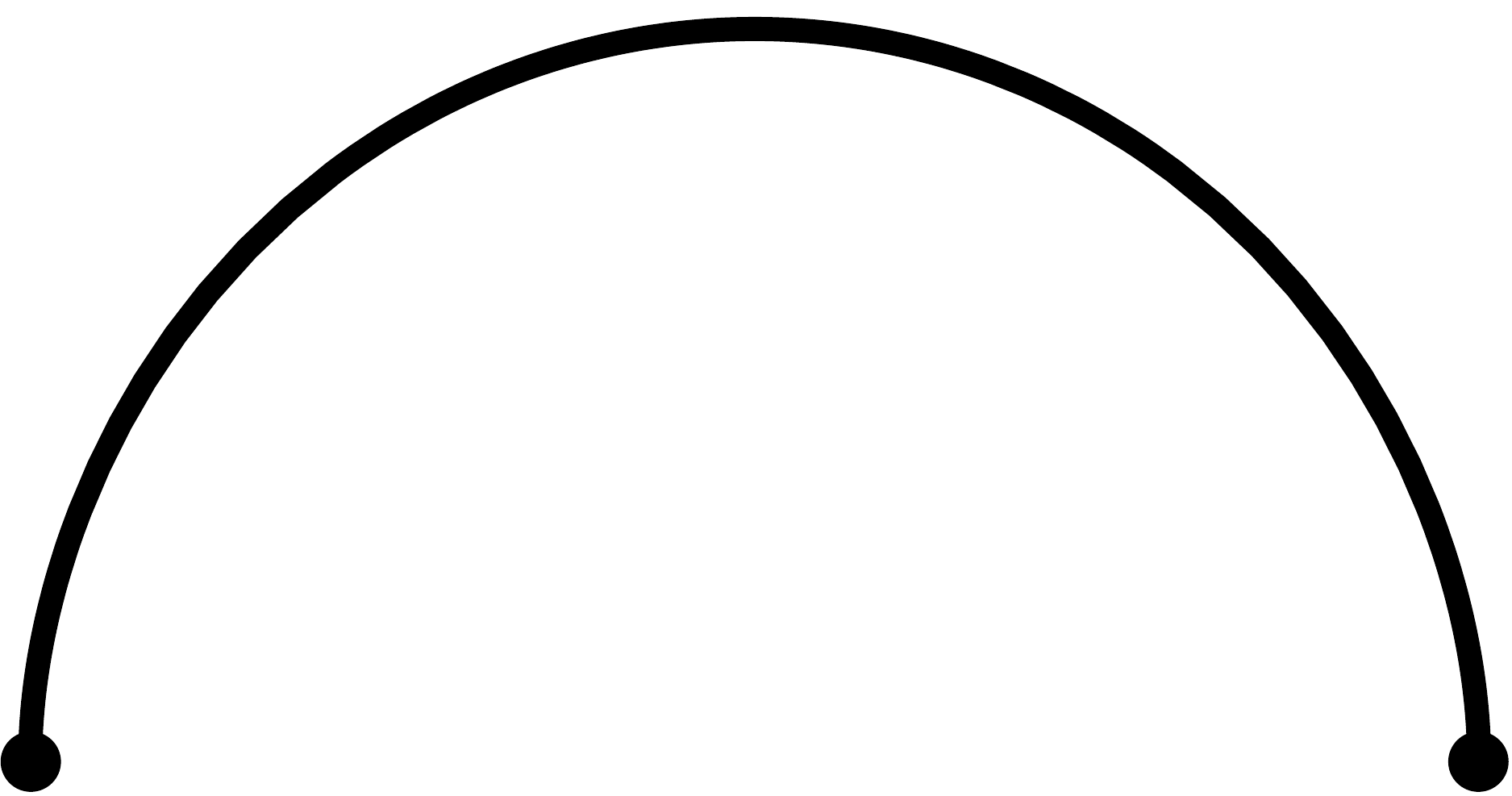}
  \endminipage\hfill
  \minipage{0.25\textwidth}
    \centering
    \includegraphics[width=0.8\linewidth]{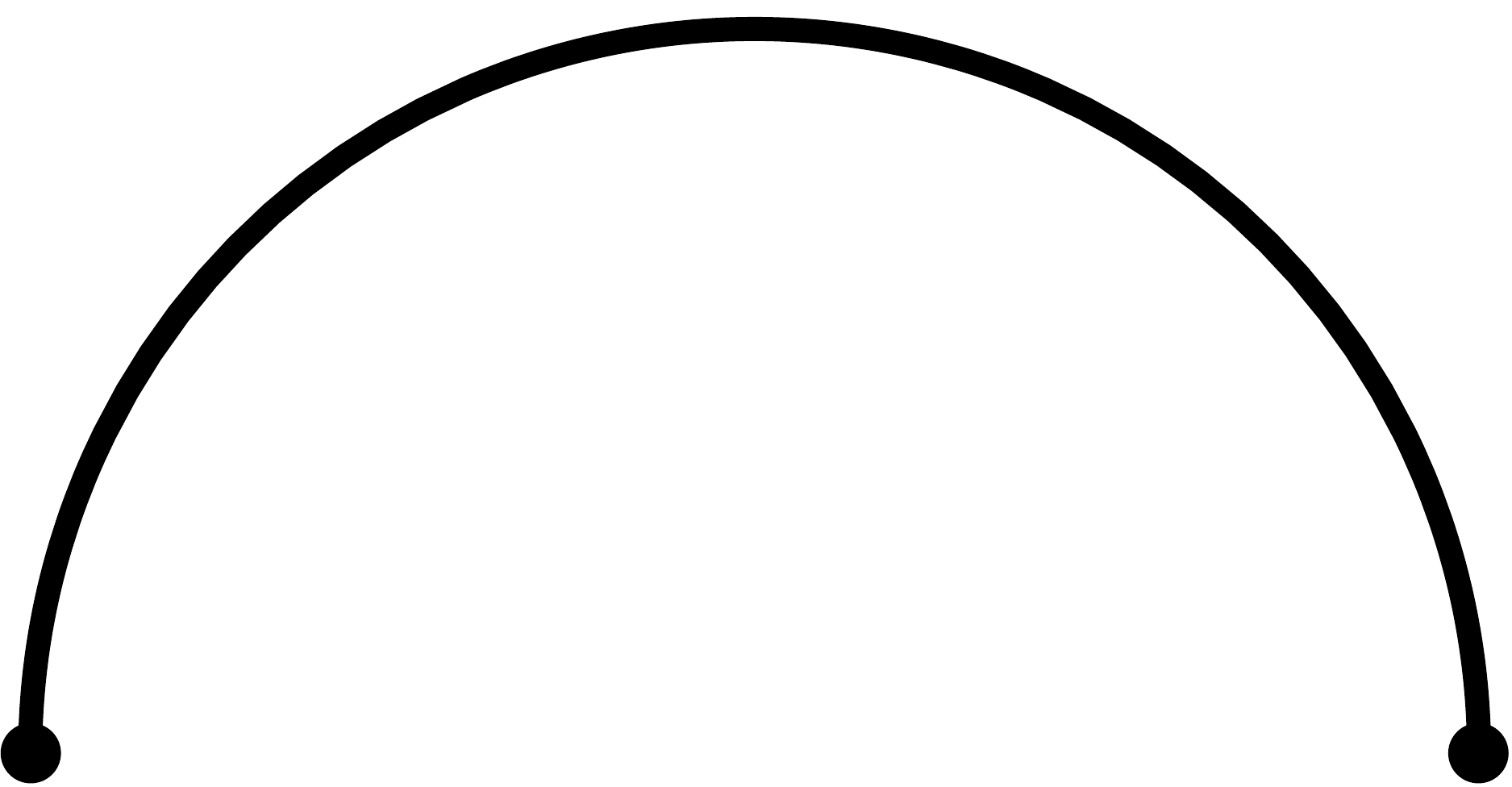}
  \endminipage\hfill
  \minipage{0.25\textwidth}
    \centering
    \includegraphics[width=0.8\linewidth]{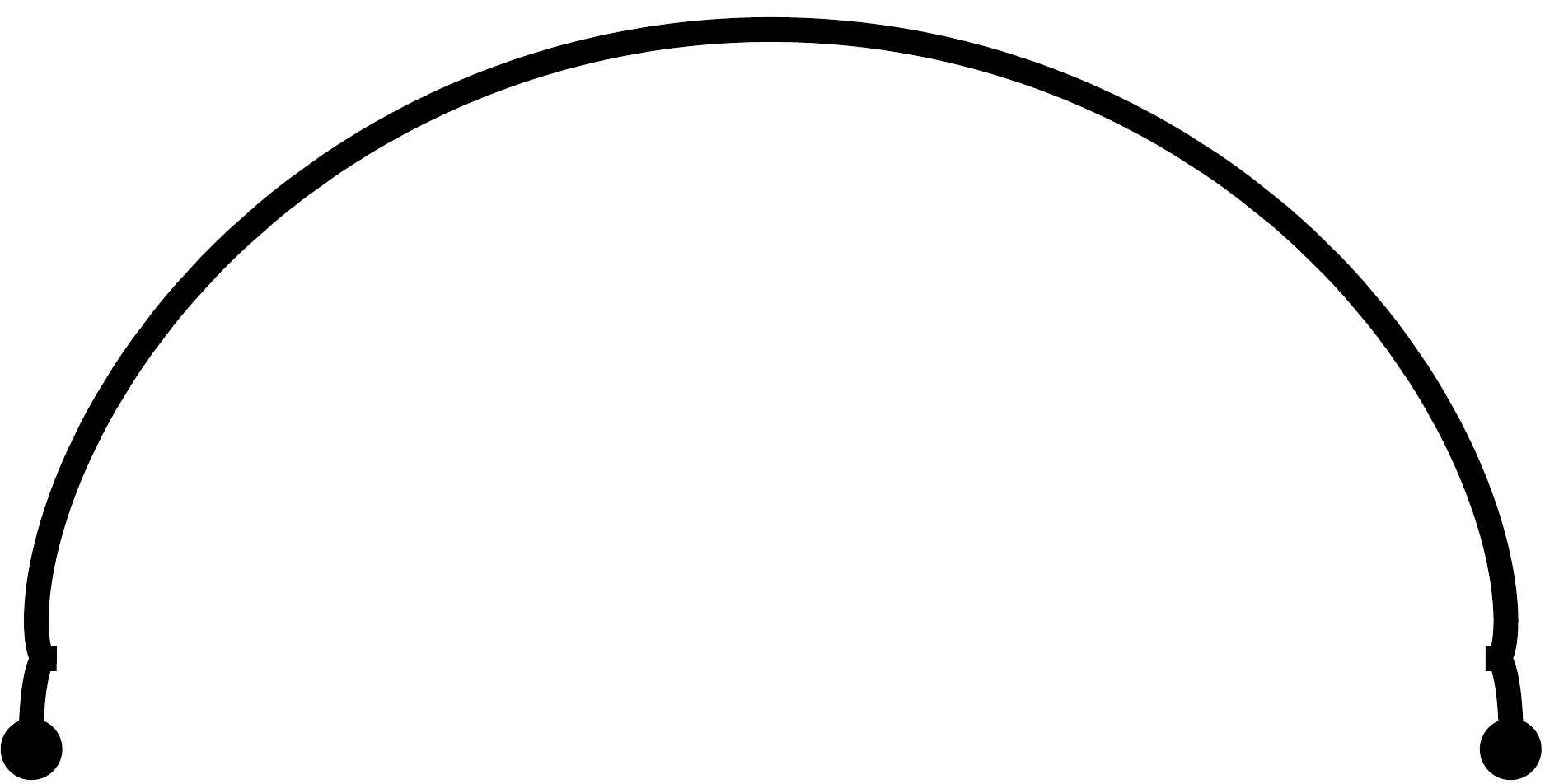}
  \endminipage\hfill
  \minipage{0.25\textwidth}
    \centering
    \includegraphics[width=0.8\linewidth]{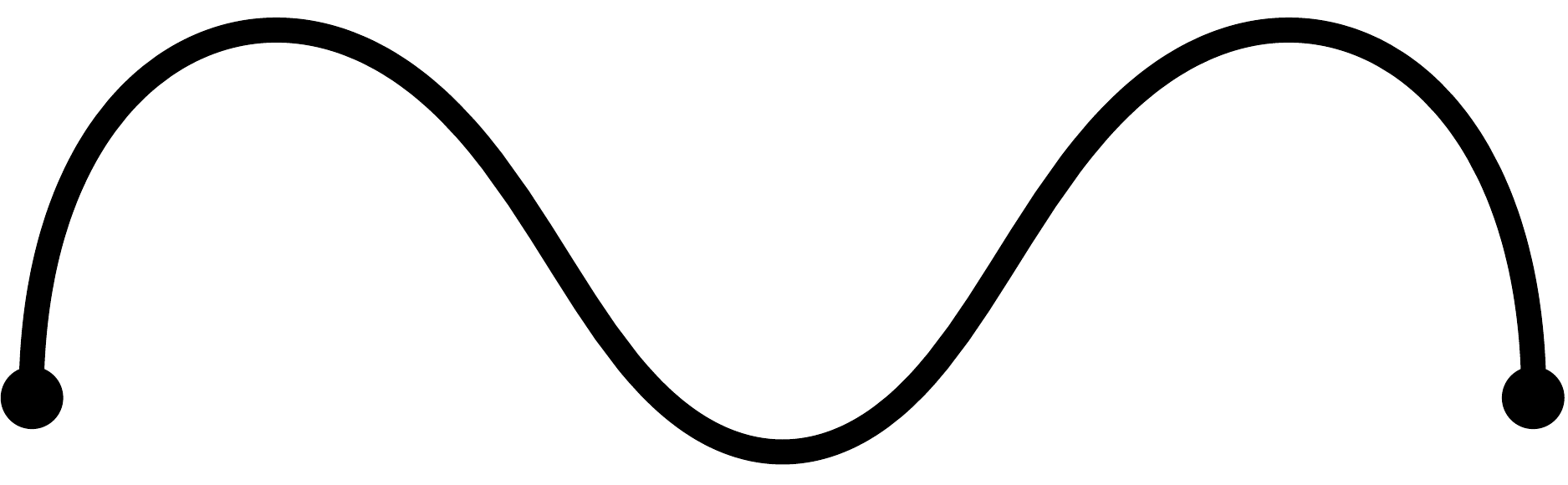}
  \endminipage\hfill \\\rule{0pt}{10ex}
  \minipage{0.25\textwidth}
    \centering
    \includegraphics[width=0.8\linewidth]{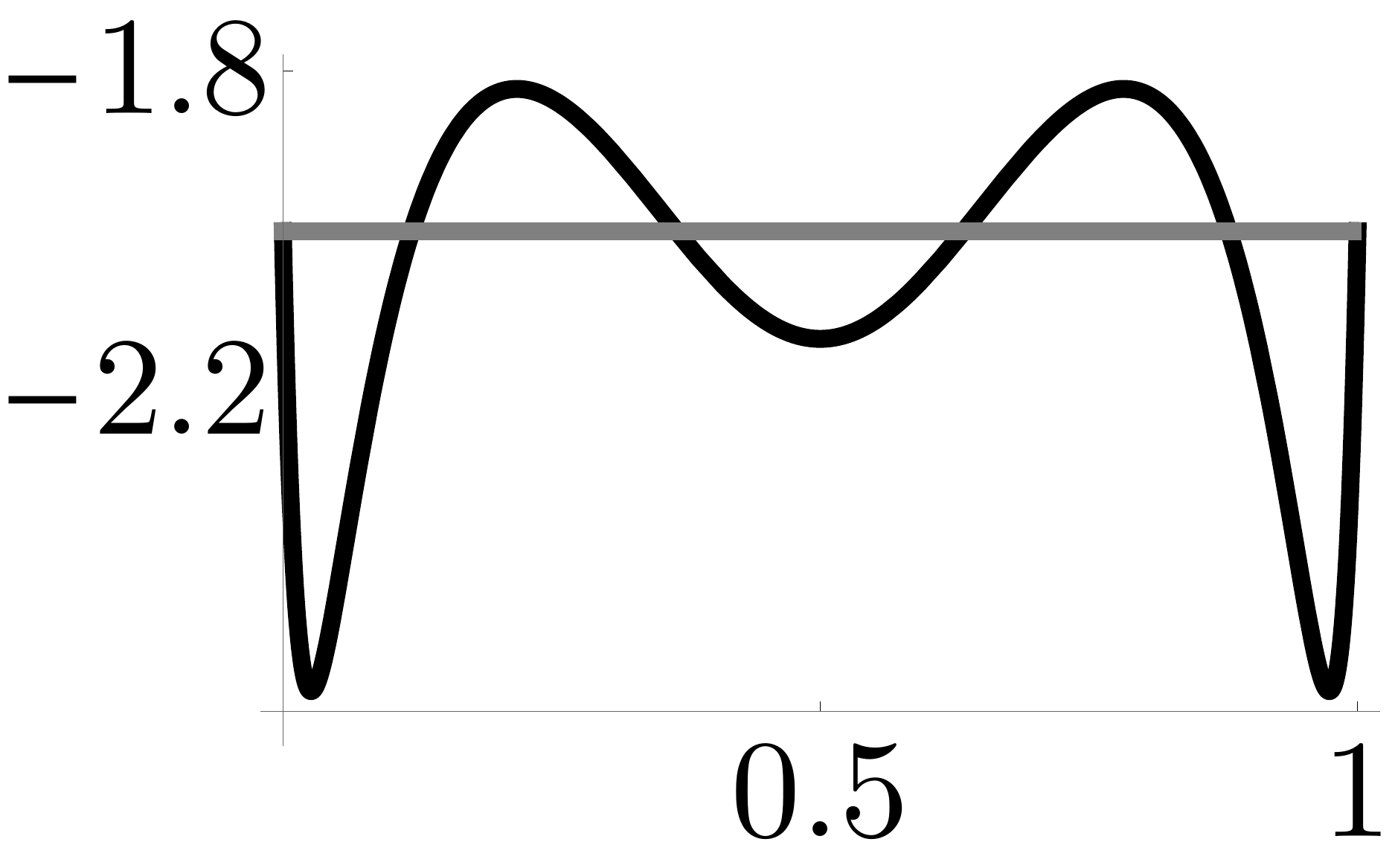}
  \endminipage\hfill
  \minipage{0.25\textwidth}
    \centering
    \includegraphics[width=0.8\linewidth]{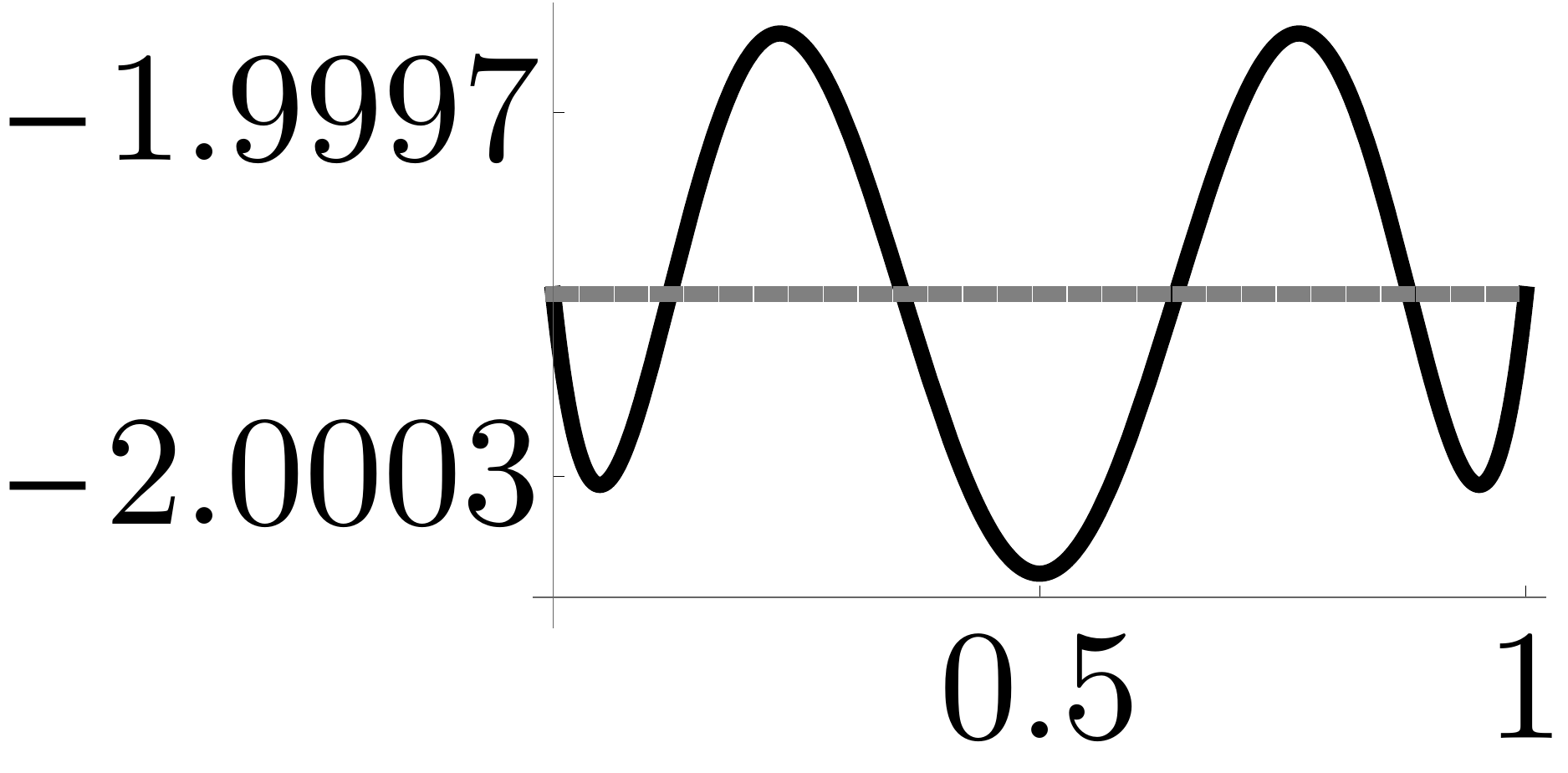}
  \endminipage\hfill
  \minipage{0.25\textwidth}
    \centering
    \includegraphics[width=0.8\linewidth]{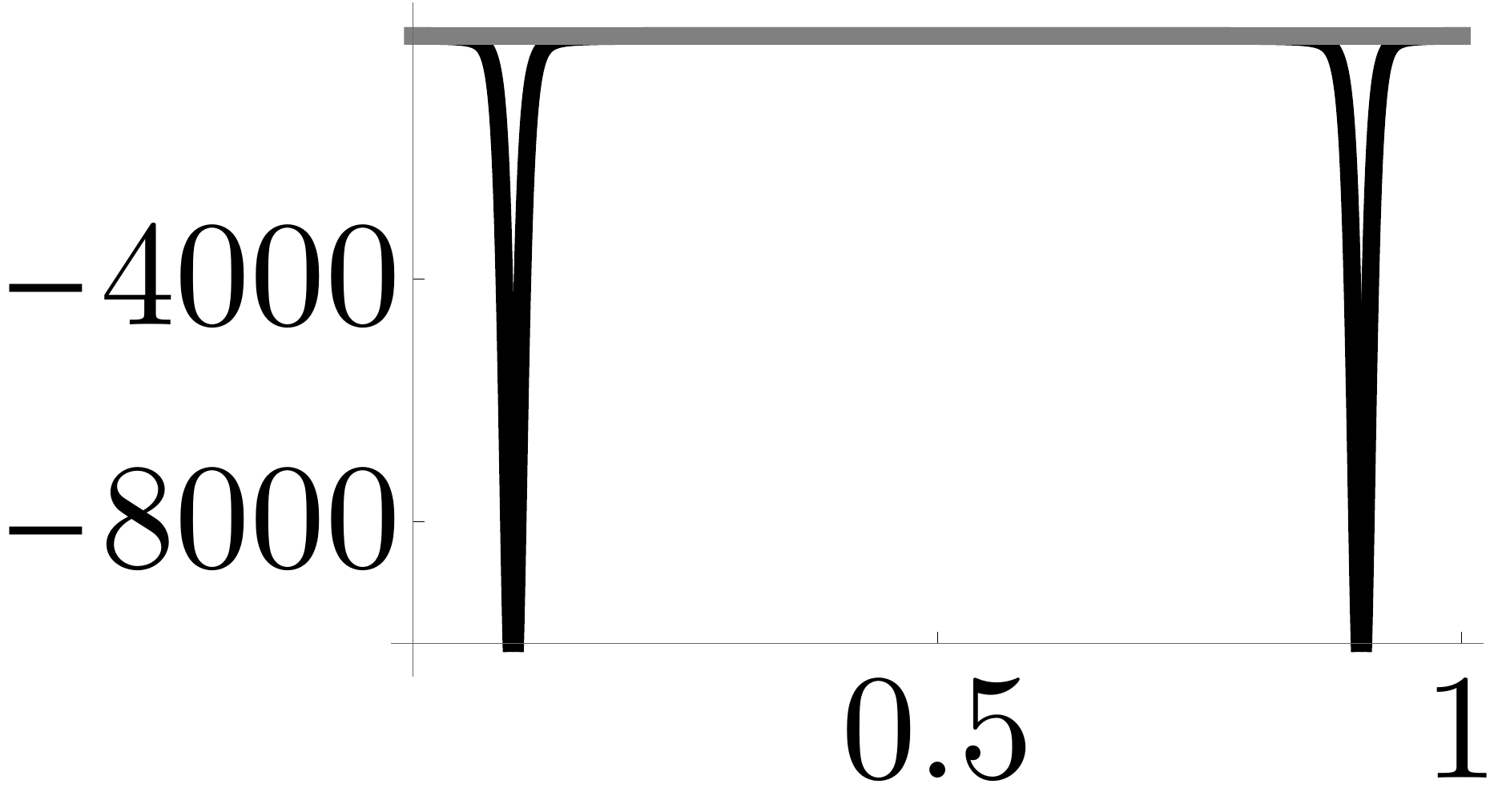}
  \endminipage\hfill
  \minipage{0.25\textwidth}
    \centering
    \includegraphics[width=0.8\linewidth]{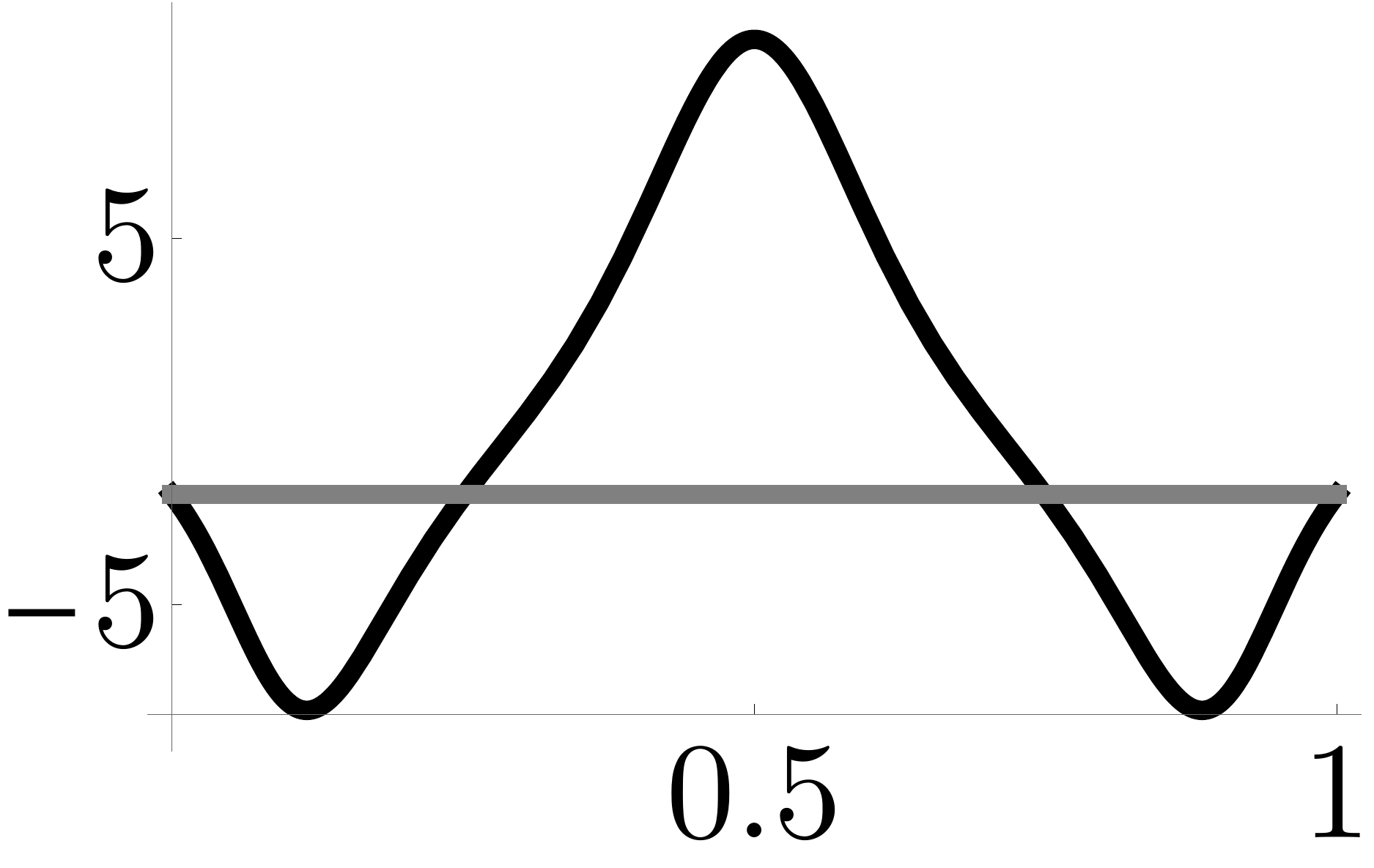}
  \endminipage\hfill
  \end{tabular}
  \caption{Plots of semicircle approximants $\bfm{p}_i$, 
  $i=1,2,3,4$, interpolating data \eqref{ex:semicircle} (top)
  together with corresponding curvature profiles (bottom). The gray
  horizontal line is the curvature of the approximated semicircle.
  Note that $\bfm{p}_3$ possesses two almost invisible tiny loops.
  \label{fig:semicircles}}
\end{figure}

In all numerical examples canonical data will be considered.
Since we know that there are always several solutions of the 
problem, we can find them numerically either by applying the 
continuation method (\cite{Allgower-90-Continuation-Methods}) 
or by using Lemma \ref{lem:four_solutions_p} which provides excellent starting values
for an iterative algorithm (such as Newton-Raphson method) 
to find positive real zeros of $p$.
One of the selection criteria \eqref{eq:rotation_index}
(for general data), \eqref{eq:curvature_error} (for general
circular arc data) or \eqref{eq:radial_distance} (for circular data with
$\alpha$ small enough) is then used to identify the most pleasant
interpolant.\\
Although we analysed the solvability of the problem only for circular data, we shall first present some
examples confirming that the interpolation method can be successfully applied for general data, too.\\
First, let us consider the data 
  \begin{equation}\label{ex:1}
    \theta_0=\pi/2,\quad \theta_1=-\pi/4,\quad 
    \kappa_0=-1,\quad \kappa_1=2,\quad 
    L=1.75.
  \end{equation}
  The problem has two solutions and they are shown in 
  Fig. \ref{fig:example1}. The one without 
  loops has the absolute rotation index 
  \eqref{eq:rotation_index} approximately $3.01$, and for the one 
  with loops we have
  $R_{abs}\approx10.43$. In this case the shape measure \eqref{eq:rotation_index} clearly
  rejects it.
\\\noindent Consider now similar data as in the previous example with parallel tangent directions
  \begin{equation}\label{ex:2}
    \theta_0=\theta_1=\pi/3,\quad\kappa_0=-1,\quad \kappa_1=2,\quad L=1.5.
  \end{equation}
  There are again two interpolants which are plotted in Fig. \ref{fig:example2}.
  In this case both of them have visually pleasant shape and its difficult to prefer one of them.
  This is also confirmed by the absolute rotation index, which is approximately $5.08$ for the left
  one and $6.14$ for the right one.
\\\noindent As the last example of general data let us consider
  \begin{equation}\label{ex:3}
    \theta_0=\pi/2,\quad \theta_1=-\pi/2,\quad \kappa_0=\kappa_1=-8,\quad L=1.21106.
  \end{equation}
  Again two formally admissible interpolants exist. They are shown in 
  Fig. \ref{fig:example3}.
  The data is actually taken from the elliptic arc $4(x-1/2)^2+16y^2=1$. It is clearly seen that the first
  approximant follows the ellipse
  (the elliptic arc is gray dashed and almost identical as the 
  interpolant), 
  while the other one takes quite different shape. The absolute rotation
  indices are $3.14$ (the analytical value of the rotation index of
  the considered elliptic arc is $\pi$) and $5.77$. The second one
  is not so much bigger since the curve is still visually pleasant. 
\\\noindent For the last two examples let us take the data from the circular arc. Suppose first that
\begin{equation}\label{ex:semicircle}
  \alpha:=\theta_0=-\theta_1=\pi/2,\quad \kappa_0=\kappa_1=-2,\quad L=\pi/2.
\end{equation}
We know from Lemma \ref{lem:four_solutions_p}  that four admissible
solutions exist. 
The error \eqref{eq:curvature_error} is taken as the selection criterion. Its corresponding values for the approximants $\bfm{p}_i$ arising
from the positive solutions $d_i$, $i=1,2,3,4$, are 
$4.2527\times 10^{-2}$, 
$8.6586\times10^{-8}$,
$2.4235\times 10^{6}$ and $34.0648$, respectiveliy.
The approximant $\bfm{p}_2$ corresponding to $d_2\approx 1.2756$ 
is clearly the most appropriate, which is confirmed also in 
Fig. \ref{fig:semicircles}.
The Hausdorff distance of the chosen interpolating curve and the circular arc is approximately
$1.2850\times10^{-5}$ and it is attained at the middle of the arc. 
Thus the constructed PH curve of degree
7 can be considered as a very accurate approximation of the semicircle
preserving the arc length.
For the $G^2$ approximation
of the whole circle just consider the spline approximant build by the constructed interpolant and its rotation.

\noindent For the last example the data from the circular arc with $\alpha>\pi/2$ will be considered. Let
\begin{equation}\label{ex:bigarc}
  \alpha:=\theta_0=-\theta_1=5\pi/6,\quad \kappa_0=\kappa_1=-1,\quad L=5\pi/3.
\end{equation}
The system of nonlinear equations \eqref{eq:curve1} and \eqref{eq:curve2} has only two 
admissible solutions. 
\begin{figure}[!htb]
  \centering
  \includegraphics[width=1\linewidth]{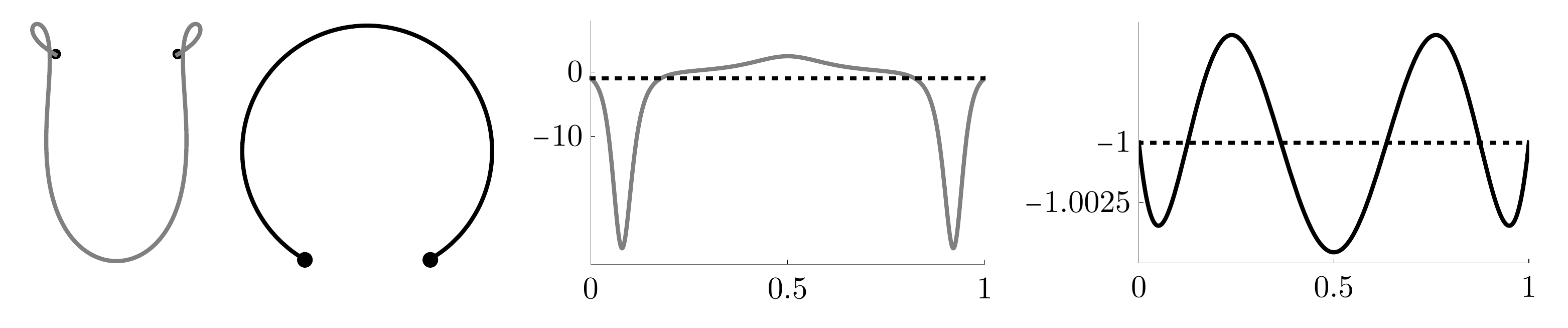}
  \caption{Plots of interpolants of data \eqref{ex:bigarc} (the first and
  the second on the left) and corresponding curvature profiles. Dotted lines
  are curvatures of the circular arc from which the data were taken.
  \label{fig:bigarc}}
\end{figure}
According to \eqref{eq:curvature_error}, the first
one is clearly rejected, since the error is 
$E_\kappa\approx61.3568$, much higher that the error of the second one
$E_\kappa\approx 9.0995\times 10^{-6}$. This is evidently confirmed also in 
Fig. \ref{fig:bigarc} where approximants together with their curvature
profiles are shown. The Hausdorff distance of the second interpolant 
and the circular arcs is
$1.6607\times 10^{-3}$, which is less than $0.2\%$ relatively to 
the radius.\\
Note that we could take even greater value of $\alpha$.
Numerical examples confirm admissible solutions of the system of
nonlinear equations for any $\alpha<\pi$, i.e., 
for the data arising from circular arcs up to almost the whole circle.\\
Let us conclude this section by numerical evaluation of the
approximation orders \eqref{eq:order2}--\eqref{eq:order34}. For the sake
of simplicity we will consider just interpolants $\bfm{p}_2$ and
$\bfm{p}_3$. They are computed for 
$\alpha_n=\pi/2^n$, $n=1,2,\dots,5$,
and the corresponding Hausdorff errors $e_n$ are determined.
From $e_n\approx c \alpha_n^r$ one easily concludes that
$r\approx \log(e_{n}/e_{n+1})/\log{2}$. Numerical results 
are collected in Tab. \ref{tab:approximation_order} and they
confirm theoretical values established in the previous section.
\renewcommand{\arraystretch}{1.5}
\begin{table}[htb]
  \begin{equation*}
    \begin{array}{|r|r|r|r|r|}\hline
       \multicolumn{1}{|c|}{\alpha}& 
       \multicolumn{1}{|c|}{d_{rad}(\bfm{p}_2;\alpha)}&
       \multicolumn{1}{|c|}{r}&
       \multicolumn{1}{|c|}{d_{rad}(\bfm{p}_3;\alpha)}&
       \multicolumn{1}{|c|}{r}\\ \hline \hline
       \pi/2 & 1.2850\times 10^{-5} & -&
       1.3865\times 10^{-2} & -\\ \hline
       \pi/4 & 6.8517\times 10^{-8} & 7.55 &
       1.3143\times 10^{-2} & 0.08\\ \hline
       \pi/8 & 4.9016\times 10^{-10} & 7.13&
       6.7687\times 10^{-3} & 0.96\\ \hline
       \pi/16 & 3.7474\times 10^{-12}& 7.03&
       3.3944\times 10^{-3} & 1.00\\ \hline
       \pi/32& 2.9119\times 10^{-14} & 7.01 &
       1.6980\times 10^{-3} & 1.00\\ \hline
    \end{array}
  \end{equation*}
  \caption{Radial distances and estimated approximation orders
  for $G^2$ interpolants $\bfm{p}_2$ (the second and the third column)
  and $\bfm{p}_3$ (the fourth and the fifth column.}
  \label{tab:approximation_order}
\end{table}
\section{Closure}\label{sec:closure}
PH curves of degree seven are promising object for interpolation 
of $G^2$ local data and preserving an arc length. Since for general data
the problem turns out to be quite complicated, a relaxation to the
circular arc data was done and a detailed analysis provided. It turned out 
that the above mentioned curves provide an excellent approximants
of circular arcs and the preserve a prescribed arc length.
An algorithm for the construction of such curves was provided.
It basically requires just solving an algebraic equation of degree six.
An asymptotic analysis reveals that the approximation order
is seven.\\
For the future work it would be nice to do some progress in studying 
the interpolation of general data. This requires some deeper analysis
of general system of nonlinear equations \eqref{eq:system}. 
Another approach to solve
the same problem would be using the PH quintic biarcs, a generalization
of cubic biarcs studied already in an early paper by \cite{Farouki-Peters-PH-cubic-biracs-1994}.

\section*{Acknowledgement}
The author was supported in part 
by the program P1-0288 and the grants J1-9104, N1-0137 and J1-3005 by 
Slovenian Research Agency.


\end{document}